\documentclass[11pt]{article}
\usepackage{mathrsfs}
\usepackage{amssymb}
\usepackage{amsmath}
\usepackage{amsbsy}
\usepackage{epsfig}
\usepackage{enumerate}
\usepackage{bm}

\RequirePackage[OT1]{fontenc}

\RequirePackage[aos]{imsart}

\newtheorem{lemma}{Lemma}[section]
\newtheorem{theorem}{Theorem}[section]

\newtheorem{assumption}{Assumption}[section]
\newtheorem{corollary}{Corollary}[section]
\newtheorem{remark}{\sc Remark}[section]

\def\proclaim#1{\par \smallskip\noindent {\bf #1}\bgroup\it\ }
\def\endproclaim{\egroup\par\smallskip}
\startlocaldefs
\endlocaldefs
\newbox\TempBox \newbox\TempBoxA

\def\pr{\textsf{P}} 
\def\ep{\textsf{E}} 
\def\Cov{\textsf{Cov}} 
\def\Var{\textsf{Var}} 

\def\Cal#1{{\mathcal #1}}
\def\bk#1{\bm #1} 

\def\underwiggle 1{
\ifmmode\setbox\TempBox=\hbox{$ 1$}\else\setbox\TempBox=\hbox{1}\fi
\setbox\TempBoxA=\hbox to \wd\TempBox{\hss\char'176\hss}
\rlap{\copy\TempBox}\smash{\lower9pt\hbox{\copy\TempBoxA}}}
\begin{document}

\begin{frontmatter}

\title{Immigrated Urn Models -- Theoretical Properties and Applications}
\runtitle{immigrated urn models}

\begin{aug}
\author{\fnms{LI-XIN} \snm{ZHANG}\thanksref{t1}\ead[label=e1]{stazlx@zju.edu.cn}},
\author{\fnms{FEIFANG} \snm{HU}\thanksref{t2}\ead[label=e2]{fh6e@virginia.edu}},
\author{\fnms{SIU HUNG}
\snm{CHEUNG}\thanksref{t3}\ead[label=e3]{shcheung@sta.cuhk.edu.hk}}
\and
\author{\fnms{WAI SUM} \snm{CHAN} \ead[label=e4]{chanws@cuhk.edu.hk}}
\address{L-X. ZHANG\\
DEPARTMENT OF MATHEMATICS$~~~~~~~~~~~~~~~$ \\
ZHEJIANG UNIVERSITY \\
HANGZHOU 310027 \\
PEOPLE'S REPUBLIC OF CHINA \\
\printead{e1}}

\address{F. HU \\
DEPARTMENT OF STATISTICS \\
UNIVERSITY OF VIRGINA\\
HALSEY HALL, CHARLOTTESVILLE~~~~~~~~~ \\
VIRGINIA 22904-4135, USA \\
\printead{e2}}

\address{S. H. CHEUNG \\
DEPARTMENT OF STATISTICS \\
THE CHINESE UNIVERSITY OF HONG KONG \\
SHATIN, N.T., HONG KONG \\
PEOPLE'S REPUBLIC OF CHINA \\
\printead{e3}}

\address{W. S. CHAN \\
DEPARTMENT OF FINANCE \\
THE CHINESE UNIVERSITY OF HONG KONG \\
SHATIN, N.T., HONG KONG \\
PEOPLE'S REPUBLIC OF CHINA \\
\printead{e4}}

\affiliation{Zhejiang University, University of Virginia and The
Chinese University of Hong Kong}

\thankstext{t1}{Research supported by  grants  from the National Natural Science Foundation of China (No. 11071214)),
 Natural Science Foundation of Zhejiang Province (No. R6100119) and  Fundamental Research Funds for the Central University (No. 2010QNA3032).}
\thankstext{t2}{Research supported by grants DMS-0349048 and DMS-0907297
from the National Science Foundation (USA).}
\thankstext{t3}{Research supported by a  grant from the Research Grants Council of the Hong Kong Special
Administrative Region (Project no. CUHK400608).}

\runauthor{L. Zhang et al.}
\end{aug}
\begin{abstract}
Urn models have been widely studied and applied in both scientific
and social science disciplines.  In clinical studies, the adoption of urn
models in treatment allocation schemes has been proved to be
beneficial to both researchers, by providing more efficient clinical
trials, and patients, by increasing the likelihood of receiving the
better treatment. In this paper, we propose a new and general class
of immigrated urn (IMU) models that incorporates the immigration
mechanism into the urn process. Theoretical properties are developed
and the advantages of the IMU models are discussed. In general, the
IMU models have smaller variabilities than the classical urn models,
yielding more powerful statistical inferences in applications.
 Illustrative
examples are presented to demonstrate the wide applicability of the
IMU models. The proposed IMU framework, including many popular classical
urn models, not only offers a unify perspective for us to comprehend
the urn process, but also enables us to generate several novel urn
models with desirable properties.

\end{abstract}



\end{frontmatter}

\section{Introduction}
\setcounter{equation}{0}

\begin{center}
\end{center}

\vspace{-0.3in}

1.1. {\em Urn Models and their applications.} Urn models have long
been considered powerful mathematical instruments in many areas,
including the physical sciences, biological sciences, social
sciences, and engineering (Johnson and Kotz, 1977; Kotz and
Balakrishnan, 1997). For example, in medical science, Knoblauch,
Neitz, and Neitz (2006) apply an urn model to study cone ratios in
human and macaque retinas. In population genetics, Hoppe (1984) and
Donnely and Kurtz (1996) employ a P\'olya-like urn model to study
Ewen's sampling distribution in neutral genetics models. Bena\"{i}m,
Schreiber, and Tarr\`{e}s (2004) also make use of a class of
generalized P\'{o}lya urn models to scrutinize evolutionary
processes.  In economics, Beggs (2005) uses the models to capture
the mechanism of reinforcement learning. In addition, numerous
examples of applications of urn models in the areas of physics,
communication theory, and computer science are  provided by
Milenkovic and Compton (2004).

In statistics, an important application  of  urn models is to
randomize treatments to patients in a clinical trial (Hu and
Rosenberger, 2006).  Consider an urn containing balls of $K$ types,
representing $K$ treatments.  Patients normally arrive sequentially,
and treatment assignment based on urn models is usually an adaptive
scheme that depends on the urn composition and previous treatment
outcomes. The urn composition is also continuously revised according
to treatment outcomes.

Early studies of urn models in statistics include the generalized
P\'{o}lya urn models (GPU) of Athreya and Karlin (1968), Wei and
Durham (1978), and Wei (1979). Another renowned variation of the
P\'{o}lya urn is the randomized P\'{o}lya urn (RPU)  proposed
by Durham, Flournoy, and Li (1998). These classic urn models have a
number of drawbacks. (i) They are usually proposed for binary
(multinomial) responses. (ii) The urn process has a predetermined
limit of urn proportions that does not have any connection with
formal optimal properties (Hu and Rosenberger, 2006). (iii) The urn
process usually has higher variability than other types of procedures
(Hu and Rosenberger, 2003) and is thus less powerful in
statistical inferences. (iv) The formulation of the asymptotic
variability is usually quite complex and it is
intricate to derive a reasonable estimate. For instance, the
asymptotic variabilities of the P\'{o}lya-urn-type models are related
to the variance of a complicated Gaussian process. In particular, for
the multi-treatment case, to derive the  variability requires
extremely complicated calculations of matrices (c.f., Smythe (1996),
Janson (2004), Bai and Hu (2005), Zhang, Hu and Cheung (2006),
Higueras et al. (2006)). (v) The models are designed mainly for the
comparison of two treatments, so there is a shortage of
methodology to handle cases with multiple treatments.

By embedding the urn process in a continuous-time birth and death
process (Ivanova et al., 2000; Ivanova and Flournoy, 2001; Ivanova,
2006), Ivanova (2003) formulates the drop-the-loser (DL) rule for a
clinical trial with two treatments. The DL rule utilizes the idea of
immigration and has been shown to yield a smaller variability among
various urn models (Hu and Rosenberger, 2003).  The DL rule is
generalized by Zhang et al. (2007) to provide more flexible urn
models.   However, these recent proposals are fragmented, offering
only a partial solution to  the aforementioned drawbacks of the
classic urn models. To supply a complete resolution, we seek to
provide a comprehensive paradigm through which one will be able to
connect existing urn models, develop useful theoretic results, and
compare merits of different classes of urn models.

\vglue 0.1in

1.2. {\em Objectives and organization of the paper.} In this paper,
we propose the IMU framework that encompasses a  wide spectrum of
urn models and incorporates the immigration process, offering a
greater flexibility in the choice of appropriate urn models in
applications. This framework includes many urn models in the
literature and provides a basis for us to derive several new urn
models, together with their desirable properties. These new urn
models are found to be capable of solving the aforementioned problems
of classic urn models.

In the literature, the asymptotic properties of urn models are
usually obtained by using Athreya and Ney's (1972) technique of
embedding the urn process in a continuous-time branching process.
However, this technique relies on the assumption that the transition
of urn composition is governed by the adding rules, which are
identical and non-random (homogeneous).  This assumption is no
longer valid for the IMU models in general due to the possibility
that the urn composition may be generated by a nonhomogeneous
immigration process. Hence, alternative mathematical approaches have
to be utilized. Another major theoretical intricacy regarding the
IMU process is that it depends on both the immigration rates and the
adding rules (refer to Section 2.1 for details).  To overcome these
mathematical difficulties, we put forward a feasible solution.
First, the IMU process is approximated by using martingales, which
can handle both immigration rates and adding rules simultaneously;
then, the IMU process is approximated by the Wiener process. Based
on the Wiener process, we will be able to obtain the asymptotic
properties of the IMU.

\vglue 0.05in To summarize, the major contributions of this paper
are as follows.
\begin{enumerate}
\item[(a)] It formulates a general framework of urn models (IMU models) that
not only encompasses most existing urn models for adaptive designs in the literature, but
also enables us to derive  new urn models with desirable properties
such as the freedom to design an urn process according to
pre-specified optimality requirements.
\item[(b)]
The paper derives asymptotic properties of the IMU models, including
strong consistency and asymptotic normality of treatment allocation
proportions. These asymptotic properties cover many existing
asymptotic properties of urn models as special cases and  form the
basis for comparisons of different IMU models.
\item[(c)] The paper proposes and discusses several  new  IMU models that are useful
in clinical trial applications.
\end{enumerate}

The  general IMU models and their asymptotic properties  are provided
in Section 2. In addition, several popular urn models that are
members of the IMU class are discussed.  In Section 3, new IMU models
are developed and their applications are given. Concluding remarks
are presented in Section 4. Finally, technical proofs are provided in
the appendix.

\vglue 0.15in
\section{The immigrated urn model}\label{section3}

\begin{center}
\end{center}

\vspace{-0.3in}

2.1. {\em The basic IMU framework.} \setcounter{equation}{0} In a
clinical trial, suppose that subjects arrive sequentially to be
randomized to one  the $K$ available treatments, and responses are
obtained immediately after treatment. An IMU model is defined as
follows. Consider an urn that contains balls of $K+1$ types. Balls of
types $1,\ldots, K$ represent treatments, and balls of type $0$ are
the immigration balls. The urn allows negative and fractional number
of balls.

Initially, there are $Z_{0,i}(\geq 0)$ balls of type $i$,
$i=0,\ldots, K$. Let $\bm Z_0=(Z_{0,0},\ldots,Z_{0,K})$ be the
initial urn composition.  Immediately before the $m$-th ($m>0$)
subject arrives to be randomized to a treatment, let the urn
composition be $\bm Z_{m-1}=(Z_{m-1,0},\ldots,Z_{m-1,K})$.  To avoid
a negative likelihood of selecting a treatment, we adopt a slight
adjustment to $\bm Z_{m-1}$ and let $Z_{m-1,i}^+=\max(0,
Z_{m-1,i})$, $i=1,\ldots, K$, and $\bm
Z_{m-1}^+=(Z_{m-1,0}^+,\ldots,Z_{m-1,K}^+)$.

To randomize the $m$-th subject, a ball is drawn at random  without
replacement. The probability of selecting a ball of type $i$ is
$Z_{m-1,i}^+/|\bm Z_{m-1}^+|$, $i=0,1,\ldots, K$. Here, $|\bm
Z_{m-1}^+|=\sum_{j=0}^K Z_{m-1,j}^+$, and $\bm Z_{m-1}^+/|\bm
Z_{m-1}^+|$ is defined to be $(0,1/K,\ldots,1/K)$ if $|\bm
Z_{m-1}^+|=0$. Hence, the balls with  negative values in $\bm
Z_{m-1}$ will have no chance of being selected unless  all
$Z_{m-1,k}^+$, $k=0,\ldots, K$ are zeros, and when $|\bm
Z_{m-1}^+|=0$ (only for the particular case where the IMU model has
no immigration ball), a treatment ball is drawn with an equal
probability of $1/K$.  Now, consider the following two
possibilities.

\begin{itemize}
\item[(a)] If the selected ball is of type $0$ (i.e, an immigration ball), no
treatment is assigned and the ball is returned to the urn.
$A_{m-1}=a_{m-1,1}+\ldots+a_{m-1,K}$ additional balls,
$a_{m-1,k}(\ge 0)$ of treatment type $k$, $k=1,\ldots, K$ are added
to the urn. Then, a ball is drawn from this updated urn again until a
treatment ball is drawn. If the immigration ball is selected $l$
times before a treatment ball is drawn, the urn composition $\bm
Z_{m-1}$ is updated to
$(Z_{m-1,0},Z_{m-1,1}+la_{m-1,1},\ldots,Z_{m-1,K}+la_{m-1,K})$ and
the $\bm Z_{m-1}^+$ is updated to
$(Z_{m-1,0},(Z_{m-1,1}+la_{m-1,1})^+,\ldots,(Z_{m-1,K}+la_{m-1,K})^+).$

\item[(b)] If a treatment ball is drawn (say, of type $k$,
$k=1,\ldots,K$), the $m$-th subject is given treatment $k$ and
the treatment outcome (response) $\xi_{m,k}$ of this subject on
treatment $k$ is observed.
 The ball is not
replaced. Instead,  $D_{m,kj}=D_{kj}(\xi_{m,k})$ balls of type $j$ are
added to the urn, $j=1,\ldots,K$. $D_{m,kj}<0$ signifies the
removal of balls.
\end{itemize}

With the IMU, the number of immigration balls remains unchanged and a
treatment ball is dropped when it is drawn. The number of treatment
balls that is added to the urn depends on: \begin{enumerate}
\item[(a)] the value of $a_{m,k}$ when an immigration ball is drawn
from the urn; and \item[(b)] the value of $D_{m,kj}$ when a
 ball of treatment type $k$ is selected.
 \end{enumerate}

Here, $a_{m,k}$s represent the immigration rates and $D_{m,kj}$s
represent the adding rules. Both $a_{m,k}$ and $D_{m,kj}$ allow
fractional values, which enable us to define a design in a flexible
manner for application. The IMU models unify many existing
urn models in the literature. Classic urn models, mainly designed
for binary responses, are members of the IMU family.  Here we list a
few popular models.

\begin{itemize}
\item[(1)] The randomized play-the-winner (RPW) rule (Wei and Durham, 1978).
When $K=2$, $Z_{0,0}=0$ or $a_{m,k}=0$ for all $m$ and $k$.  Further,
$D_{m,kk}=2$ if the response of the $m$-th subject on
treatment $k$ is a success, and
 $D_{m,kk}=D_{m,kj}=1$ ($j \neq k$) otherwise.
\item[(2)] Generalized P\'{o}lya urn models (Athreya and Karlin, 1968,
also called the generalized Friedman's urn). When $a_{m,k}=0$, we
obtain the GPU models if one chooses the adding rule $D_{m,kj}$ as in
Section 4.1 in Hu and Rosenberger (2006). If $D_{m,kj}$ is
non-homogeneous, we obtain the non-homogeneous GPU models discussed
by Bai and Hu (1999, 2005).
\item[(3)] The birth and death urn (BDU) (Ivanova et al., 2000).
Suppose that $a_{m,k}\equiv 1$, $D_{m,kj}=0$ for $j\ne k$.  In addition,
$D_{m,kk}=2$ if the response of the $m$-th subject on
treatment $k$ is a success, and $D_{m,kk}=0$ otherwise.
When $K=2$, we obtain the birth and death urn (BDU) (Ivanova et al., 2000).
When $K>2$, we obtain generalized birth and death urn (BDU) for $K$ treatments.
\item[(4)] The Drop-the-loser (DL) rule (Ivanova, 2003).
Suppose that $a_{m,k}\equiv 1$, $D_{m,kj}=0$ for $j\ne k$.  In addition,
$D_{m,kk}=1$ if the response of the $m$-th subject on
treatment $k$ is a success, and
$D_{m,kk}=0$ otherwise.
When $K=2$, we obtain the DL rule (Ivanova, 2003).
When $K>2$, we obtain DL rule for $K$ treatments.
\item[(5)] The generalized drop-the-loser (GDL) rule (Zhang et al., 2007).
Suppose that $a_{m,k}=a_k$ (does not depend on $m$) are constants and
$D_{m,kj}=0$ for $j\ne k$. When $K=2$, we obtain the GDL rule. When
$K>2$, we obtain GDL rules for $K$ treatments.
\item[(6)] Sequential estimated urn (SEU) models (Zhang, Hu and Cheung, 2006).
When $Z_{0,0}=0$ or $a_{m,k}=0$ for all $m$ and $k$, and $D_{m,kj}$ depends on estimation, we
obtain the SEU models proposed by Zhang, Hu and Cheung (2006) and the urn models in Bai, Hu and Shen (2002).
\end{itemize}
In general, we can select suitable $a_{m,k}$ and $D_{m,kj}$ to obtain
the desirable IMU model for both binary and continuous responses (see
examples in Section 3).

In clinical trials, let $N_{n,k}$ be the number of subjects who have
been assigned to treatment $k$, $k=1,\ldots, K$. Denote $\bm
N_n=(N_{n,1},\ldots,N_{n,K})$. In clinical studies, the proportions
$N_{n,k}/n$, $k=1,\ldots, K$ of patients being assigned to various
treatments are  useful statistics. In fact, for urn model
applications, there are several important statistics, including:

\begin{itemize}
\item[(a)] the urn proportion $Z_{n,k}/\sum_{k=1}^K Z_{n,k}$;
\item[(b)] the allocation proportion $N_{n,k}/n$; and
\item[(c)] the estimation of the unknown parameters in the model.
\end{itemize}

It is worthwhile  noting that both $a_{m,k}$ and $D_{m,kj}$ depend on
$m$. This allows both the immigration rates and the adding rules to
be expressed as functions of all previous responses thus far in the
clinical trial. Then, we are able to construct  desirable IMU
models that can be used to suit pre-specified  allocation proportion
targets. To reiterate, as both $a_{m,k}$ and $D_{m,kj}$ depend on
$m$, it is impossible to use Atheya and Ney's (1972) technique of
embedding the urn process in a continuous-time branching process.

It is also worth noting that  Hoppe's urn (Hoppe, 1984) and its
extensions (see for example Donnely and Kurtz, 1996) are not members
of the IMU models. For Hoppe's urn, the number of ball types is
increasing and random, but for an IMU model the number of ball types
is fixed ($K + 1$).

\smallskip

2.2. {\em Notation and assumptions.} Before the discussion of major
asymptotic results regarding the IMU models, we  introduce some basic
notation and the necessary assumptions. Suppose that $\xi_{m,k}$
($k=1,...,K$, $m=1,2,3,...$) is the random variable representing the
response of the $m$-th subject on  treatment $k$. In practice, we
only observe one $\xi_{m,k}$ for each $m$. Without loss of
generality, we assume that the unknown parameter $\theta_k$ is the
mean of the outcome $\xi_{m,k}$ and  take the sample mean as its
estimate. Write $\bm \xi_m=(\xi_{m,1},\ldots,\xi_{m,K})$.  For the
adding rules, let $\bm D_m=(D_{m,kj}; k,j=1,\ldots,K)$, $\bm
D_m^{(k)}=(D_{m,k1},\ldots,D_{m,kK})$, $k=1,\ldots, K$,  and $\bm
H_m=(h_{kj}(m))=\ep \bm D_m$. Let $\widehat{\theta}_{m-1,k}$ be the
sample mean of the responses
\begin{equation}\label{eqestimate}
\widehat{\theta}_{m-1,k}=\frac{c_1+S_{m-1,k}}{c_2+ N_{m-1,k} },
\end{equation}
where $S_{m-1,k}$ is the sum  of the responses on treatment $k$ of
all the previous  $m-1$  subjects. Here, $c_1,~c_2>0$ are used to
avoid the nonsense case of $0/0$.  These two constants  play a
minor role, only  in the earlier stages of the clinical trial when
accumulated observations of the treatments are still very small. In
general, many estimators, such as the MLE, can be written in
the form of (\ref{eqestimate})
 with $S_{m-1,k}$ being replaced by a sum of functions
of the responses plus a negligible remainder (see Hu and Zhang
(2004a) for detail discussion).

As discussed in Section 2.1, the immigration rate $a_{m,k}$ plays an
important role in the IMU models. Its significance will be
illustrated in the later part of this section when the theoretical
properties of the IMU models are being reviewed. In clinical trials,
optimal allocation proportions usually depend on the unknown
parameters $\bm\theta$ (See Rosenberger, et al., 2001 and Tymofyeyev,
Rosenberger and Hu, 2007). To achieve these
proportions, one can select the immigration rates $a_{m,k}$ as
functions of $\bm\theta$.  In practice, as $\bm\theta$ is unknown,
one can use $a_{m-1,k}=a_k(\widehat{\bk \theta}_{m-1})$ as the
immigration rates. The guidelines for the selection of the function
$a_k$ will be given in Section 3. In most applications, the adding
rules $\bm D_m=(D_{m,kj}; k,j=1,\ldots,K)$ normally depend on the
response $\bk \xi_{m}$, similar to those in the GPU models.
Hence, we need the following assumptions.

\begin{assumption}\label{Asp0}
Functions $a_k(\cdot)>0$ are continuous and twice differentiable at
$\bm\theta$.
\end{assumption}

\begin{assumption}\label{Asp2} $\{(\xi_{m,k}, D_{m,k1},\ldots, D_{m,kK});m\ge 1\}$, $k=1,\ldots,
K$, are $K$   sequences of i.i.d. random variables with
$\sup_m\ep|D_{m,kj}|^{2+\delta}<\infty$, and
$\sup_m\ep|\xi_{m,k}|^{2+\delta}<\infty$ for some $0<\delta\le 2$,
$k=1,\ldots, K$. Hence, let $\bm H_m=\bm H$, which does not
depend on $m$. Further assume that $D_{m,kk}\ge -C$ for some
$C$, $k=1,\ldots, K$, and also $D_{m,kj}\ge 0$ for $k\ne j$.
\end{assumption}

The continuity of $a_k(\cdot)$ in Assumption \ref{Asp0} is needed to
show that $0<\min_{m,k}a_{m,k}\le \max_{m,k} a_{m,k}<\infty$ as given
in Lemma \ref{lem6}. The differentiability of the function is
required for the Taylor expansion. The moment condition in Assumption
2.2 is useful for applying the limit theorems and the approximation
of related martingales. Finally, the lower bound of $D_{m,kj}$
implies that when a ball is drawn, the maximum number of balls of
that treatment type which can be removed is $C+1$. This condition is
used to derive the lower bound of $Z_{n,k}$, as given in Lemma
\ref{lem4}.

\vglue 0.1in 2.3. {\em Main asymptotic results.}  We now discuss the
asymptotic properties related to urn proportions and model parameter
estimators.  Asymptotic results are classified into one of following
three possible cases, according to the expectation of the adding
rules.

\begin{enumerate}
\item $\bm H\bm 1^{\prime}<\bm 1^{\prime}$ where $\bm
1=(1,\ldots,1)$.  Hence $\sum_{j=1}^Kh_{kj}<1$ for all
 $k=1,\ldots, K$.  The urn composition is mainly updated by the
 immigration balls because, on average, the number
of added balls in each step according to the outcome of a treatment
is less than the number of dropped balls, which is 1. The derivation
of asymptotic results for this case is of the utmost importance and
plays a crucial role in this paper.
\item
$\bm H\bm 1^{\prime}>\bm 1^{\prime}$. The total number of balls in
the urn gradually increases to infinity.  Hence, the probability of
drawing an immigration ball drops to zero. For this case, we will
prove that the IMU model is asymptotically equivalent to the generalized
P\'{o}lya urn model without immigration (refer to Theorem
\ref{theoremGPU}).
\item
$\bm H\bm 1^{\prime}=\bm 1^{\prime}$.  This is the borderline
case in which both the treatment balls and the immigration ball
retain their roles in the urn updating process.
\end{enumerate}

These three cases lead to very different asymptotic results. Let us
first consider the case of $\bm H\bm 1^{\prime}>\bm 1^{\prime}$. The
following theorem  ensures that the IMU model behaves asymptotically,
the same as the generalized P\'{o}lya urn model,  when $\bm H\bm
1^{\prime}>\bm 1^{\prime}$. The proof is given in the appendix.
Based on this theorem, we can obtain the asymptotic properties,
including the strong consistency, asymptotic normality and Gaussian
approximation, of the generalized P\'{o}lya urn model as discussed
by Janson (2004), Bai and Hu (2005), Zhang, Hu and
Cheung (2006), Zhang and Hu (2009), among others.

\begin{theorem} \label{theoremGPU} Suppose that Assumption  \ref{Asp2}
is satisfied, $\bm H\bm 1^{\prime}=\gamma\bm 1^{\prime}$ with
$\gamma>1$, and $0\le a_{m,k}\le C m^{1/2-\delta_0}$ for some
$\delta_0>0$ and all $m,k$. Let $\bm v=(v_1,\cdots, v_K)$  be the
left eigenvalue vector of $\bm H$ that corresponds to the largest
eigenvalue $\gamma$ and satisfies $v_1+\cdots+v_K=1$, and denote
$\widetilde{\bm H}=\frac{\bm H-\bm I}{\gamma-1}-\bm 1^{\prime}\bm
v$. Further, let $\lambda_2,\cdots,\lambda_K$ be the other $K-1$
eigenvalues of $\bm H$ and
$\lambda=\max\{Re(\lambda_2),\cdots,Re(\lambda_K)\}$. Assume that
$\lambda-1<(\gamma-1)/2$. Then, there exist two independent standard
$K$-dimensional Wiener processes $\bm B_{t1}$ and $\bm B_{t2}$ such
that
$$ (N_{n,1},\cdots, N_{n,K})-n\bm v=\bm G_{n1}+\frac{1}{\gamma-1}\int_0^t \frac{\bm G_{x2}}{x}dx (\bm I-\bm 1^{\prime}\bm v)+o(n^{1/2-\epsilon}) \;a.s., $$
$$ (Z_{n,1},\cdots,Z_{n,K})-(\gamma-1)n\bm v=(\gamma-1)\bm G_{n1}\widetilde{\bm H}+\bm G_{n2} +o(n^{1/2-\epsilon}) \;a.s.,$$
for some $\epsilon>0$, where $\bm G_{ti}$ is the solution of the equation
$$ \bm G_{ti}=\bm B_{ti}\bm \Lambda_i^{1/2}+\int_0^t \frac{\bm G_{xi}}{x} dx \widetilde{\bm H}, $$
with $\bm \Lambda_1=diag(\bm v)-\bm v^{\prime}\bm v$ and $\bm\Lambda_2=\sum_{k=1}^K v_k\Var\{\bm D_1^{(k)}\}$.
In particular,
$$  \frac{Z_{n,0}}{Z_{n,0}+\cdots+Z_{n,K}}\to 0\; a.s., \;\; \frac{Z_{n,k}}{Z_{n,0}+\cdots+Z_{n,K}}\to v_k\; a.s., \;\; \frac{N_{n,k}}{n}\to v_k\; a.s.,$$
 $k=1,\cdots, K$, and
$$ n^{1/2}\left(\frac{Z_{n,1}}{(\gamma-1)n}-v_1,\cdots,\frac{Z_{n,K}}{(\gamma-1)n} -v_K\right)\overset{\mathscr{D}}\to N\left(\bm 0,\bm\Gamma_1\right), $$
$$ n^{1/2}\left(\frac{N_{n,1}}{n}-v_1,\cdots,\frac{N_{n,K}}{n}- v_K\right)\overset{\mathscr{D}}\to N\left(\bm 0,\bm\Gamma_2\right). $$
Here, the variance-covariance matrices $\bm \Gamma_1$ and $\bm
\Gamma_2$ can be  specified in line with Bai and Hu (2005) and Zhang and Hu
(2009) with $\frac{\bm D_m-\bm I}{\gamma-1}$ and $\frac{\bm H-\bm
I}{\gamma-1}$ replacing $\bm D_m$ and $\bm H$, respectively. For
details, one can refer to Proposition 3.4 of Zhang and Hu (2009).
\end{theorem}

Now we  consider the case in which $\bm H\bm 1^{\prime}<\bm 1^{\prime}$.
Different from the case when $\bm H\bm 1^{\prime}>\bm 1^{\prime}$ in
which the urn proportion and the sample allocation proportion have
the same limit, the urn proportion may not have a limit in this case.
For the immigration rates, write $a_k=a_k(\bm \theta)$.  Let $\bm
a=(a_1,\ldots, a_K)$, $\bm u=\bm a(\bm I-\bm H)^{-1}$, $s=\bm a(\bm
I-\bm H)^{-1}\bm 1^{\prime}=\sum_{k=1}^Ku_k$ and $\bm v=\bm u/s$.
Further, denote $\bm \Sigma_k=\Var\{\bm D_1^{(k)}\}$, $\bm
\Sigma_{11}=\sum_{k=1}^K v_k \bm \Sigma_k$, $\bm
\Sigma_{12}=(\Cov\{D_{1,kj},\xi_k\}; j,k=1,\ldots, K)$, $\bm
\Sigma_{22}=diag(\Var\{\xi_{1,1}\},\ldots, \Var\{\xi_{1,K}\})$,  and
\begin{equation}\label{eqLambda}
 \bm \Lambda=
 \begin{pmatrix} \bm \Lambda_{11} & \bm \Lambda_{12}  \\
                             \bm \Lambda_{12}^{\prime} & \bm \Lambda_{22} \end{pmatrix}
 =\begin{pmatrix} \bm \Sigma_{11} & \bm \Sigma_{12}\; diag(\bm v) \\
                             diag(\bm v)\; \bm \Sigma_{12}^{\prime} & \bm \Sigma_{22}\; diag(\bm
                             v)\end{pmatrix}.
\end{equation}

\begin{theorem}\label{th3} Suppose that Assumptions \ref{Asp0}-\ref{Asp2} are
satisfied, $\bm H\bm 1^{\prime}<\bm 1^{\prime}$ and $Z_{0,0}>0$.  Then
$$ Z_{n,k}=o(n^{1/2-\epsilon})\;\; a.s., \;\; k=1,\ldots, K$$
for some $\epsilon>0$, and, one can define a $2K$-dimensional
Wiener processes $(\bm W(t), \bm B(t))$  such that
\begin{equation}\label{eqth3.1}
\Var\{(\bm W(t),
\bm B(t))\}=t\bm \Lambda
\end{equation}
and
\begin{align}\label{eqth3.2}
\bm N_n-n\bm v=\bm W(n)\bm A
 +\int_0^{n}\frac{\bm B(x)}{x}dx\;diag\Big(\frac{1}{\bm v} \Big)
 \frac{\partial \bm v(\bm\theta)}{\partial\bm\theta}+o(n^{1/2-\epsilon})\;\; a.s.
\end{align}
for some $\epsilon>0$, where $\bm A=(\bm I-\bm H)^{-1}(\bm I-\bm 1^{\prime}\bm v)$,
$$ \bm v=\bm v(\bm\theta)=\frac{a(\bm \theta)(\bm I-\bm H)^{-1}}{a(\bm \theta)(\bm I-\bm H)^{-1}\bm 1^{\prime}}
\;\; \text{ and } \;\;\frac{\partial \bm
v(\bm\theta)}{\partial\bm\theta}=\left(\frac{\partial
v_k(\bm\theta)}{\partial \theta_j};j,k=1,\ldots,K\right).$$ Here,
$1/\bm v=(1/v_1,\ldots, 1/v_K)$.
\end{theorem}

\begin{remark} Note that $h_{ij}\ge 0$ for $i\ne j$. The existence of
$(\bm I-\bm H)^{-1}$ is implied by the assumption that $\bm H \bm
1^{\prime}<\bm 1^{\prime}$.
This assumption  can be replaced by a more general assumption
in which there is a vector $\bm e= (e_1,\ldots,e_K)$ such that $\bm
H\bm e^{\prime}<\bm e^{\prime}$  and $e_i>0$, $i=1,\ldots, K$.
\end{remark}

Based on Theorem \ref{th3}, we can see that the urn composition
$(\sqrt{n})^{-1}Z_{n,k}$ converges to $0$ almost surely. It is because
when $\bm H\bm 1^{\prime}<\bm 1^{\prime}$, there will be a
net loss of balls from the urn on average if a treatment ball is
drawn. The proof of Theorem \ref{th3} is given in the appendix. The
consistency  and asymptotic normality of $\bm N_n$ can be
derived by using Equation (\ref{eqth3.2}) as follows.
\begin{corollary}\label{normal}Under the assumptions in Theorem \ref{th3},
\begin{equation}\label{eqnormal1.1}
 \bm N_n - n\bm v =O(\sqrt{n\log\log n})\; a.s. \;  \text{ and } \;
\sqrt{n}\big(\frac{\bm N_n}{n}-\bm v\big)\overset{\mathscr{D}}\to N(\bm
0,\bm\Sigma),
\end{equation}
where
$\bm\Sigma=\bm\Sigma_{D}+2\bm\Sigma_{\xi}+\bm\Sigma_{D\xi}+\bm\Sigma_{D\xi}^{\prime}$,
and
\begin{align*}
 &\bm\Sigma_D=\bm A^{\prime}\bm\Sigma_{11}
\bm A,\quad \bm\Sigma_{D\xi}=\bm A^{\prime} \bm\Sigma_{12}  \frac{\partial \bm v(\bm\theta)}{\partial\bm\theta}\\
&\bm\Sigma_{\xi}=\Big(\frac{\partial \bm v(\bm\theta
)}{\partial\bm\theta}\Big)^{\prime} diag\Big(\frac{\Var\{\xi_{1,1}\}}{v_1},\ldots,\frac{\Var\{\xi_{1,K}\}}{ v_K}\Big)\frac{\partial
\bm v(\bm\theta)}{\partial\bm\theta}.
\end{align*}
In particular, if $\bm D_m \equiv const$, then
$$ \sqrt{n}\big(\frac{\bm N_n}{n}-\bm v\big)\overset{\mathscr D}\to N(\bm 0,2\bm\Sigma_{\xi});
$$
and if $a_{m,k}\equiv a_k$, $k=1,\ldots,K$, do not depend on the
estimates, then
$$ \sqrt{n}\big(\frac{\bm N_n}{n}-\bm v\big)\overset{\mathscr D}\to N(\bm 0,\bm\Sigma_{D}).
$$
\end{corollary}
{\bf Proof}. $~$ Note that $(\bm W(n),\int_0^n\frac{\bm B(x)}{x}dx)$ is
a centered Gaussian vector with
$$\bm W(n)=O(\sqrt{n\log\log n})\;\; a.s., $$
$$\int_0^n\frac{\bm B(x)}{x}dx=O(1)+
\int_e^n\frac{O(\sqrt{x\log\log x})}{x}dx=O(\sqrt{n\log\log n})\;\; a.s.,$$
$$ \Var\{\bm W(n)\}=n \bm\Sigma_{11},
$$
\begin{align*}
\Var\left\{\int_0^n\frac{\bm B(x)}{x}dx\right\}=
\bm\Sigma_{22} diag(\bm v)\int_0^n\int_0^n\frac{x\wedge y}{xy}dxdy =2n
\;\bm\Sigma_{22} diag(\bm v)
\end{align*}
and
$$\Cov\left\{\bm W(n),\int_0^n\frac{\bm B(x)}{x}dx\right\}=
\bm\Sigma_{12}diag(\bm v)\int_0^n\frac{x\wedge n}{x}dx=n\;\bm\Sigma_{12}diag(\bm v).
$$
(\ref{eqnormal1.1}) follows from (\ref{eqth3.2}) immediately.
$\Box$

\begin{remark}
In practice, the responses in clinical trials are frequently not
available immediately before the treatment allocation of the next
subject (delayed response). The parameters can be estimated and the
urn can be updated only by using all available observed responses.
In the delayed response situation,  let $\mu_k(m,l)$ be
the probability that
 the response of the $m$-th
subject on treatment $k$ occurs after at least another $l$ subjects
arrive. If $\mu_k(m,l)\le C l^{-\gamma}$ for some $\gamma\ge 2$, then we
can show that the total sum of unobserved outcomes up to the $n$-th
assignment
 is with a high order of $\sqrt{n}$ and thus the conclusion in
 Theorem \ref{th3}  remains true.
It has been shown that the delay mechanism does not effect the
asymptotic properties for many response-adaptive  designs if the
delay decays with a power rate (c.f., Bai, Hu, and Rosenberger,
2002; Hu and Zhang, 2004b; Zhang et al., 2007).
\end{remark}

In many IMU models (such as, special cases (3), (4) and (5) in
Section 2), the additional rule, $\bm D_m$, is a diagonal matrix
($D_{m,kj}=0$, $j\ne k$). For this special case, we have the
following corollary that helps us to obtain the asymptotical limits
and covariance matrix of $\bm N_n$  easily.

\begin{corollary}\label{gdl}Suppose that Assumptions \ref{Asp0}-\ref{Asp2} are
satisfied,  $D_{m,kj}=0$ for $j\ne k$, and $h_k=1-\ep D_{1,kk}>0$.
Write $\bm h=(h_1,\ldots, h_K)$,
$$v_k(\bm\theta,\bm h)=\frac{a_k(\bm\theta)/h_k}{\sum_{j=1}^Ka_j(\bm\theta)/h_j}\;\;
k=1,\ldots, K,
$$
 $\bm v=\bm v(\bm\theta,\bm h)=(v_1(\bm\theta,\bm h),\ldots,v_K(\bm\theta,\bm h))$, and
\begin{align*}
 \frac{\partial
\bm v(\bm\theta,\bm h)}{\partial\bm\theta} =&\left(\frac{\partial
v_k(\bm\theta,\bm h)}{\partial \theta_j};j,k=1,\ldots,K\right),
\quad \\
 \frac{\partial \bm v(\bm\theta,\bm h)}{\partial\bm h}
 =&\left(\frac{\partial v_k(\bm\theta,\bm h)}{\partial
h_j};j,k=1,\ldots,K\right).
\end{align*}
Then,
\begin{equation}\label{eqCor1.2}
\frac{\bm N_n}{n}\to \bm v \quad a.s. \;  \text{ and } \;
\sqrt{n}\big(\frac{\bm N_n}{n}-\bm v\big)\overset{\mathscr D}\to N(\bm
0,\bm\Sigma),
\end{equation}
where
$\bm\Sigma=\bm\Sigma_{D}+2\bm\Sigma_{\xi}+\bm\Sigma_{D\xi}+\bm\Sigma_{D\xi}^{\prime}$,
\begin{align*}
 &\bm\Sigma_D=\Big(\frac{\partial \bm v(\bm\theta,\bm h)}{\partial\bm h}\Big)^{\prime} diag\Big(\frac{\sigma_{D
1}^2}{v_1},\ldots,\frac{\sigma_{D K}^2}{ v_K}\Big)\frac{\partial
\bm v(\bm\theta,\bm h)}{\partial\bm h},\\
&\bm\Sigma_{\xi}=\Big(\frac{\partial \bm v(\bm\theta,\bm
h)}{\partial\bm\theta}\Big)^{\prime} diag\Big(\frac{\sigma_{\xi
1}^2}{v_1},\ldots,\frac{\sigma_{\xi K}^2}{ v_K}\Big)\frac{\partial
\bm v(\bm\theta,\bm h)}{\partial\bm\theta},
\\
&\bm\Sigma_{D\xi}=-\Big(\frac{\partial \bm v(\bm\theta,\bm
h)}{\partial\bm h}\Big)^{\prime} diag\Big(\frac{\sigma_{D\xi
1}}{v_1},\ldots,\frac{\sigma_{D\xi K}}{ v_K}\Big)\frac{\partial \bm
v(\bm\theta,\bm h)}{\partial\bm\theta},
\end{align*}
and $\sigma_{Dk}^2=\Var\{D_{1,kk}\}$, $\sigma_{\xi
k}^2=\Var\{\xi_{1,k}\}$, $\sigma_{\xi Dk}= \Cov\{D_{1,kk},
\xi_{k,1}\}$, $k=1,2,\ldots,K$.
\end{corollary}
{\bf Proof.}
It is easy to check that
$$\bm \Sigma_{11}=diag(\sigma_{D1}^2v_1,\ldots,\sigma_{DK}^2v_K),~
\bm \Sigma_{12}=
  diag(\sigma_{\xi D1},\ldots,\sigma_{\xi DK}),$$
$$  \bm \Sigma_{22}=
diag(\sigma_{\xi 1}^2,\ldots,\sigma_{\xi  K}^2),~ \bm A=diag(1/\bm
h)(\bm I-\bm 1^{\prime}\bm v),$$ and $ \partial \bm v(\bm\theta,\bm
h)/ \partial\bm h=-diag(\bm v)\bm A$. Then, the results follow from
Corollary \ref{normal} directly. $\Box$

\smallskip To improve statistical efficiency, a suitable response
adaptive randomization procedure should be adopted because of
variability (Hu and Rosenberger, 2003). Hu, Rosenberger and Zhang
(2006) studied the variability of a randomization procedure that
targets any given allocation proportion. They obtained a lower bound
of the variability. For a large class of the IMU models in this
paper, the lower bound of the variability is attained. When the
variance of IMU model attains the lower bound, we can use the
Cram\'er-Rao formula to compute the variance. In general, we have the
following theorem.

\begin{theorem}\label{theorem4.1}
If each $D_{m,kj}$ is a linear function of
 a random $\eta_{m,k}$, $j=1,\ldots, K$, where $\eta_{m,k}$ may be a function of $\xi_{m,k}$ and for each $k$, $\eta_{m,k}$, $m=1,2,\ldots, $ are i.i.d.
 random variables with finite variances, then we have
 \begin{align} \bm \Sigma_D=&\left(\frac{\partial \bm v}{\partial \bm d}\right)^{\prime}
diag\left(\frac{\Var\{\eta_{1,1}\}}{v_1},\ldots,\frac{\Var\{\eta_{1,K}\}}{v_K}\right)
\frac{\partial \bm v}{\partial \bm d} \label{eqlowboundV1}\\
\bm \Sigma_{D\xi}=&\left(\frac{\partial \bm v}{\partial \bm
d}\right)^{\prime}
diag\left(\frac{\Cov\{\eta_{1,1},\xi_{1,1}\}}{v_1},\ldots,\frac{\Cov\{\eta_{1,K},\xi_{1,K}\}}{v_K}\right)
\frac{\partial \bm v}{\partial \bm \theta}, \label{eqlowboundV2}
\end{align}
where $\bm d=(d_1,\ldots,d_K)=(\ep\eta_{1,1},\ldots,\ep\eta_{1,K})$.
Further, if $\bm a(\cdot)=const$ and $\Var\{\eta_{1,k}\}$ is the
inverse of the Fisher information of $d_k$, then the asymptotic
variance-covariance matrix  of $\bm N_n/\sqrt{n}$ attains the
following lower bound,
\begin{equation}\label{lowerb}
\left(\frac{\partial \bm v}{\partial \bm d}\right)^{\prime}
diag\left((v_1I_1)^{-1},\ldots,(v_K I_K)^{-1}\right)
\left(\frac{\partial \bm v}{\partial \bm d}\right),
\end{equation}
where $I_k$ is the Fisher information function of parameter $d_k$.
\end{theorem}

 {\bf Proof.} If we write $\bm D_1^{(k)}=\bm
\alpha_k+\bm\beta_k\eta_{1,k}$ and $\bm K=\begin{pmatrix} \bm \beta_1
\\ \cdots \\\bm \beta_K\end{pmatrix}$, then
\begin{align*} \bm\Lambda_{11}=&\sum_{k=1}^K v_k\Var\{\eta_{1,k}\} \bm\beta_k^{\prime}\bm\beta_k \\
=&(diag(\bm v)\bm K)^{\prime}
diag\left(\frac{\Var\{\eta_{1,1}\}}{v_1},\ldots,\frac{\Var\{\eta_{1,K}\}}{v_K}\right)diag(\bm
v)\bm K
\end{align*}
and
\begin{align*} \bm\Sigma_{12}= (diag(\bm v)\bm K)^{\prime}
diag\left(\frac{\Cov\{\eta_{1,1},\xi_{1,1}\}}{v_1},\ldots,\frac{\Cov\{\eta_{1,K},\xi_{1,K}\}}{v_K}\right).
\end{align*}
However, $\partial \bm H/\partial d_k=diag(\bm 1_k)\bm K$, where
$\bm 1_k$ has zero elements except the $k$-th one which is $1$.
In addition,
$$ \frac{\partial (\bm I-\bm H)^{-1}}
{\partial d_k}=(\bm I-\bm H)^{-1}\frac{\partial \bm H}{\partial
d_k}(\bm I-\bm H)^{-1} =(\bm I-\bm H)^{-1} diag(\bm 1_k)\bm K (\bm
I-\bm H)^{-1}.$$ It follows that
\begin{align*}
\frac{\partial\bm v}{\partial d_k}=&\frac{\partial \bm a(\bm I-\bm
H)^{-1}/\partial d_k}{\bm a(\bm I-\bm H)^{-1}\bm 1^{\prime}}
-\frac{\partial\bm a (\bm I-\bm H)^{-1}/\partial d_k}{(\bm a(\bm
I-\bm H)^{-1}\bm 1^{\prime})^2}
\bm 1^{\prime} \bm a(\bm I-\bm H)^{-1} \\
= & \bm v diag(\bm 1_k)\bm K(\bk I-\bm H)^{-1}(\bm I-\bm
1^{\prime}\bm v)=\bm v diag(\bm 1_k)\bm K\bm A,
\end{align*}
i.e., $\partial \bm v/\partial \bm d=diag(\bm v)\bm K\bm A$.  Hence, (\ref{eqlowboundV1}) and (\ref{eqlowboundV2}) are proved by
Corollary \ref{normal}. $\Box$

\vglue 0.1in Corollary \ref{gdl} and Theorem \ref{theorem4.1} are
useful for deriving the asymptotic variance. We will illustrate this
idea by introducing several interesting examples in the next
section.

\begin{remark}
In Theorem \ref{theorem4.1}, for simplicity of notation we assume
that the parameter  $d_k$ is a one-dimensional  parameter that
corresponds to treatment $k$. The theorem is still valid if
reformulated using a vector parameter $d_k$, without extra
 assumptions.
\end{remark}

\smallskip

Finally, we consider the case when $\bm H\bm 1^{\prime}=\bm
1^{\prime}$. The following theorem, with proof given in the appendix,
can be used to yield the consistency property of the allocation
proportion. However, it is still unknown  whether $\bm N_n$ is
asymptotically normal.
\begin{theorem} \label{th4} Suppose that Assumptions \ref{Asp0} and \ref{Asp2} are
satisfied, and $\bm H\bm 1^{\prime}=\bm 1^{\prime}$, $Z_{0,0}>0$.
Suppose further that $1$ is a single eigenvalue of $\bm H$. Then
$$ \bm N_n-n\bm v=O(\sqrt{n\log\log n}) \;\; a.s.\;\; \text{ and } \;\; \bm N_n-n\bm v=O_P(\sqrt{n}), $$
where $\bm v$ is the left eigenvalue vector of $\bm H$ that
corresponds to the eigenvalue $1$ and satisfies $v_1+\cdots+v_K=1$.
\end{theorem}

These theorems and corollaries are related to the sample
allocation proportion $\bm N_n/n$. Regarding the estimator
$\hat{\bm\theta}_n$, we have the following theorem.
\begin{theorem}\label{ththeta}
Suppose that the assumptions in Theorem \ref{theoremGPU} or
\ref{th3} or \ref{th4} are satisfied. We have
\begin{equation}\label{eqthetaesti}
\sqrt{n}\left(\hat{\bm\theta}_n-\bm\theta\right) \to N(\bm
0,\bm\Sigma_{\bm\theta}),
\end{equation}
where
$$\bm\Sigma_{\bm\theta}=diag\Big(\frac{\Var\{\xi_{1,1}\}}{v_1},\ldots,\frac{\Var\{\xi_{1,K}\}}{ v_K}\Big). $$
\end{theorem}

Note that $\bm N_n/n\to \bm v$ a.s. according to
Theorem \ref{theoremGPU} or \ref{th3} or \ref{th4}, so the proof of this
Theorem is the same as that of Lemma 1 of Hu, Rosenberger and Zhang
(2006) and is thus omitted here.

\section{Examples and Applications}
\setcounter{equation}{0}

In this section, we apply the general asymptotic results in Section 3
to selected IMU models for illustrative purposes. In Section 2.1, we
listed several classic families of urn models as special cases of
IMU models. We can  apply directly the theoretical results in Section
3 to these special cases and obtain their asymptotic properties for
both $K=2$ (available in the literature) and for the general value of $K
\geq 3$. In this section, we focus on the generation of new families
of urn models from the IMU framework and discuss their corresponding
properties. Several illustrative examples are given. First, we
consider continuous-type responses that
 are frequently  encountered in
clinical studies, even though there has been a lack of related
studies in the literature.

\vglue 0.1in  \noindent {\it Example 1: Two treatments with
continuous responses}. Suppose that  $\xi_{m,1}$ ($m=1,2,3,...$) are i.i.d.
random variables from $N(\mu_1,\sigma_1^2)$ and  $\xi_{m,2}$
($m=1,2,3,...$) are i.i.d. random variables from
$N(\mu_2,\sigma_2^2)$.  Without the loss of generality, assume that
the smaller the value of the response, the better the treatment. We
now introduce four  IMU models.
\begin{itemize}
\item[(1.A)]
 Let $a_{m,k} \equiv 1$, $D_{m,kj}=0$ for $j\ne k$.  Let $C$
be a constant such that
$D_{m,kk}=1$ if the response of the $m$-th subject on
treatment $k$, $\xi_{m,k}$, is less than $C$, and
 $D_{m,kk}=0$ otherwise.
\item[(1.B)]
Suppose that there are two critical values $C_1<C_2$ and if it is
very desirable to have the value of  the response fall between $C_1$
and $C_2$, then the following IMU model is appropriate. Take $a_{m,k}
\equiv 1$, $D_{m,kj}=0$ for $j\ne k$.  Further, let
$D_{m,kk}=1$ if $\xi_{m,k} < C_1$,
 $D_{m,kk}=0$ if $\xi_{m,k} > C_2$ and else
 $D_{m,kk}=1/2$.
\item[(1.C)] If the power of statistical inferences is an important concern, the Neyman allocation
$\sigma_1/(\sigma_1+\sigma_2)$ can be adopted to maximize the power
of testing.  Then, consider the following IMU model. Let $a_{m,k}=
\widehat{\sigma}_k$, $D_{m,kj}=0$ for all $j,k$. Here,
$\widehat{\sigma}_k^2$ is the current sample variance of the
responses on treatment $k$, $k=1,2$, and can be used as estimates in
the Neyman allocation rule.
\item[(1.D)] If the aim is to lower the proportion of subjects being assigned to the inferior treatments
 for ethical reasons, the allocation target
$\sqrt{\mu_2} \sigma_1 /(\sqrt{\mu_2} \sigma_1+\sqrt{\mu_1}
\sigma_2)$ where $\mu_1,\mu_2>0$ (Zhang and Rosenberger (2006)) is an
option. Let $a_{m,1}= \sqrt{\widehat{\mu}_2} \widehat{\sigma}_1$,
$a_{m,2}= \sqrt{\widehat{\mu}_1} \widehat{\sigma}_2$, $D_{m,kj}=0$
for all $j,k$. Here, $\widehat{\mu}_k$, $\widehat{\sigma}_k^2$ are the
current sample mean and sample variance of the responses on treatment
$k$ respectively, $k=1,2$. To avoid the situation of
$\widehat{\mu}_k\le 0$, simply replace $\widehat{\mu}_2$ by $1/m$
when such an occasion arises.
\end{itemize}

\noindent Designs (1.A) and (1.B) cover a wide spectrum of potential
applications.  Note that Design (1.A) is equivalent to the DL rule
for binary response if the critical value $C$ is used to classify
responses into two categories.  Designs (1.C) and (1.D) incorporate
pre-specified objectives of a clinical trial, depending on whether
the objective is to increase the testing power (as in (1.C)), or
reduce the number of patients being assigned to the inferior
treatments (as in (1.D)). Further, it would not difficult to
generalize these four designs to studies with $K > 2$ treatments.

The asymptotic properties of the four designs can be obtained
using Theorem 2.2. For illustrative purposes, we discuss
asymptotic normalities for Designs (1.C).
It is easy to verify that
\begin{align*}
 \widehat{\sigma}_k^2=:\widehat{\sigma}_{m,k}^2=&\frac{1}{N_{m,k}}\sum_{j=1}^m X_{j,k}(\xi_{j,k}-\mu_k)^2
 -\big(\widehat{\mu}_k-\mu_k\big)^2 \\
 =&\frac{1}{N_{m,k}}\sum_{j=1}^m X_{j,k}(\xi_{j,k}-\mu_k)^2
 +O\big(\frac{\log\log N_{m,k}}{N_{m,k}}\Big)\;\; a.s.
\end{align*}
By Corollary \ref{normal},
$$ \frac{N_{n,1}}{n}\to v_1\; a.s. \text{ and } n^{1/2}\left(
\frac{N_{n,1}}{n}- v_1\right)\overset{\mathscr{D}}\to N(0,\sigma^2), $$
where $v_1=\sigma_1/(\sigma_1+\sigma_2)$, and $\sigma^2$ equals to
$$ 2 \left(\frac{\partial v_1}{\partial (\sigma_1^2)},\frac{\partial v_1}
{\partial (\sigma_2^2)}\right)
diag\left(\frac{\Var\{(\xi_{1,1}-\mu_1)^2\}}{v_1},\frac{\Var\{(\xi_{1,2}-\mu_2)^2\}}
{1-v_1}\right)\left(\frac{\partial v_1}{\partial
(\sigma_1^2)},\frac{\partial v_1} {\partial
(\sigma_2^2)}\right)^{\prime}. $$ After simplification, we have
$\sigma^2=\sigma_1\sigma_2/(\sigma_1+\sigma_2)^2.$ One can also use
Theorem \ref{ththeta} to derive the asymptotic distribution of the
estimators of the unknown parameters. For example, in Design (1.C),
$\sqrt{n}(\widehat{\sigma}_{n,k}^2-\sigma^2)\overset{\mathscr{D}}\to
N(0,2\sigma_k^4/v_k)$.

\vglue 0.2in \noindent {\it Example 2: Modified DL (MDL) rule}. We
propose the MDL rule, which is a modification of the DL rule.  The
procedure is similar to the DL rule in that when a treatment ball is
drawn, this ball is replaced only when the response is a success.
However, when an immigration ball is drawn, instead of adding an
equal number of treatment balls to the urn, we add $C\widehat{p}_k$
($C>0$) balls of type $k$, $k=1,\ldots,K$, where $\widehat{p}_k$ is
the current estimate of the successful probability $p_k$ of treatment
$k$, and $C$ is a constant. With this model, more balls are
immigrated to treatments with  higher success rates, and
subsequently, the limit proportions will be higher for better
treatments.

Regarding the asymptotic variance, it is straightforward to show that
$\bm a=(p_1C,\ldots,p_KC)$ and $\bm H=diag(p_1,\ldots,p_K)$. The
conditions in Corollary \ref{gdl} are satisfied for all  cases with
$0<p_k<1$ and $k=1,\ldots, K$. Hence, the limit proportions are $
v_k=(p_k/q_k)/({\sum_{j=1}^Kp_j/q_j}), \; k=1,\ldots, K. $ The
asymptotic  variance-covariance can be derived by the formulae in
Corollary \ref{gdl},  in which $\bm\theta=(p_1,\ldots,p_K)$, $\bm
h=(q_1,\ldots,q_K)$, and $\sigma_{Dk}^2=\sigma_{\xi k}^2=\sigma_{D\xi
k}=p_kq_k$, $k=1,\ldots, K$. For the two-treatment case,
$$ \frac{N_{n,1}}{n}\to v_1=\frac{p_1/q_1}{p_1/q_1+p_2/q_2}\; a.s.\;\;\text{ and }\;
\; \sqrt{n}( N_{n,1}/n- v_1)\overset{\mathscr D}\to N(0,\sigma^2),
$$ where $
\sigma^2=q_1q_2[p_1^2(1+q_2^2)+p_2^2(1+q_1^2)]/(p_2q_1+p_1q_2)^3. $
When the success probabilities $p_1$ and $p_2$ are both high, the
variability $\sigma^2$ is close to the lower bound $q_1q_2(p_1^2
+p_2^2 )/(p_2q_1+p_1q_2)^3$.

\smallskip
Unlike the generalized P\'{o}lya urn models without immigration in
which the asymptotic normality holds only when a very strict
condition on eigenvalues of a generating matrix is satisfied (c.f.,
Bai and Hu, 2005; Janson, 2004; Zhang, Hu, and Cheung, 2006),  the
MDL rule allows asymptotic
normality for all cases  with $0<p_k<1$, $k=1,\ldots, K$.

\smallskip

In most IMU models, the adding rule $\bm D_m$ is a diagonal matrix.
Here we give an example for the two-treatment case with dichotomous
responses in which the adding rule $\bm D_m$ is not a diagonal
matrix.

\vglue 0.1in \noindent {\it Example 3: Two treatments with
dichotomous responses}. Consider the two-treatment case with
dichotomous responses, success or failure. Let $p_k$ be the success
probability of  treatment $k$ and $q_k=1-p_k$, $k=1,2$. We consider
an  immigrated urn in which $a_{m,1}=a_{m,2}\equiv 1$ and
$$ \bm D_m=\begin{pmatrix} \beta \xi_{m,1} &
\alpha (1-\xi_{m,1}) \\ \alpha (1-\xi_{m,2}) & \beta \xi_{m,2}
\end{pmatrix}, $$
where $\xi_{m,k}=1$ if the outcome of the $m$-th subject on
treatment $k$ is a success, and $0$  otherwise, $k=1,2$, $\alpha\ge
0$. In this design, the draw  of an immigration ball generates a
ball of each treatment type;  when a treatment type ball is dropped,
$\beta$ balls of the same treatment type are added if the outcome is
a success and $\alpha$ balls of the alternate treatment type are
added if the outcome is a failure. Hence,
$$\bm H= \begin{pmatrix} \beta p_1 & \alpha q_1 \\ \alpha q_2 & \beta p_2 \end{pmatrix}. $$
Based on Theorems 2.1-2.5 of Section 2, we can derive the asymptotic
properties for the three cases: (i) $\bm H\bm 1^{\prime}> \bm
1^{\prime}$; (ii) $\bm H\bm 1^{\prime}< \bm 1^{\prime}$; and (iii)
$\bm H\bm 1^{\prime}= \bm 1^{\prime}$. The technical details are omitted
here. Nevertheless, it is worth  noting that different choices of
$\alpha$ and $\beta$ generate various members of the IMU family.

\begin{remark}
The GDL rule of Zhang et al. (2007) is a member of the IMU class
with $D_{m,kj}=0$, $j\ne k$. In practice, the values of $\theta_k,
k=1,...,K$ are unknown and have to be estimated by sample
statistics.  The derivation of the asymptotic distributions of the
treatment proportions $N_{n,k}$ is usually difficult and is not
included by Zhang et al. (2007) if the estimates of $\theta_k,
k=1,...,K$ are used. However, by applying Corollary \ref{gdl}, one
can obtain the asymptotic properties of $N_{n,k}$ directly.

For example, if the optimal proportion
$v_1=\sqrt{p_1}/(\sqrt{p_1}+\sqrt{p_2})$ is used for comparing two
treatments, we can select an IMU model with  $\bm D_m\equiv 0$,
$a_{m,k}=C \sqrt{\widehat{p}_k}$, where $\widehat{p}_k$ is the
current estimate of the successful probability $p_k$ of treatment
$k$, and $C$ is a constant, $k=1,2$. By Corollary \ref{gdl}, we have
$$\sqrt{n}(N_{n,1}/n-v_1)\overset{\mathscr D}\to N(0,\sigma^2), \mbox{ where } \sigma^2=\frac{1}{2(\sqrt{p_1}+\sqrt{p_2})^3}
\left(\frac{p_2q_1}{\sqrt{p_1}}+\frac{p_1q_2}{\sqrt{p_2}}\right).$$
Zhang et al. (2006) proposed the use of a GPU without immigration to
target this proportion (c.f., their Example 2). The corresponding
asymptotic variance is
$$\frac{\sqrt{p_1p_2}}{(\sqrt{p_1}+\sqrt{p_2})^2}+\frac{3}{2(\sqrt{p_1}+\sqrt{p_2})^3}\left(\frac{p_2q_1}
{\sqrt{p_1}}+\frac{p_1q_2}{\sqrt{p_2}}\right),$$ which is at least
triple the variance of this IMU model.
\end{remark}

The IMU models, such as those given in the foregoing examples, can be
applied in clinical trials. We discuss the applications in three
possible directions.
\begin{itemize}
\item[(i)]
There are numerous applications of urn models in clinical trials. One
can apply the proposed IMU models with multiple objectives, such as
ethical concerns and design efficiency. For instance, Tamura et al.
(1994) discussed the application of the RPW rule, a member of the IMU
family, to study the treatment of out-patients suffering from
depressive disorder.  Later, in a simulation study (using the same
data), Bhattacharya (2008) showed that the DL rule, another member of
the IMU family, has a smaller variability and yields higher power
than the RPW rule.  One can apply the asymptotics of the IMU model
given in this paper to compare various urn allocation methods instead
of using only the simulation results  given by Bhattacharya (2008).
\item[(ii)] Urn models are also frequently  employed in clinical studies to
promote balance (see Matthews et al. (2010) and the references therein).
In such circumstances, IMU models should be considered as useful
candidates. The introduction of the immigration urn will
significantly improve these allocation schemes, mainly in relation to the
variability of the urn proportions.  Furthermore, asymptotic
distributions of IMU models can be derived, leading to a more
comprehensive understanding of these urn processes.
\item[(iii)] For comparing $K$ treatments, Tymofyeyev, Rosenberger and Hu (2007), Zhu and Hu (2009)
obtained optimal allocation proportions for both binary and
continuous responses. The IMU models are suitable choices due to
their low variability and flexibility in targeting these optimal
allocation proportions.
\end{itemize}

\section{Conclusions}\label{section5}

In this paper, we have proposed a general class of urn models that
incorporates immigration.  The IMU framework unifies many existing
classes of urn models and provides  crucial linkages among these
models to enable us to have a more comprehensive understanding of
different urn processes and their important properties. Further, this
framework facilitates the generation of new urn models with desirable
properties. Asymptotic properties of the IMU models, with widely
satisfied conditions, are given in Section 2. These important
 results serve to connect existing asymptotic results
about urn models. More importantly, the asymptotic normality formula
in this article can be employed to evaluate and compare different urn
models in terms of the distributions of treatment allocation
proportions. Under very mild conditions, the suggested IMU models
always yield relatively smaller asymptotic variances. In many cases,
the asymptotic variance attains the lower bound. Thus, the IMU models
 have smaller variabilities than the corresponding generalized
P\'olya urn models.

In clinical trials,  responses may not be available immediately
after the patients have been treated.  However, there are no
logistical difficulties in incorporating delayed responses into the IMU
framework. One can update the urn when responses become available. A
moderate delay in response (see Hu and Zhang, 2004b) will not affect
the asymptotic properties of the IMU. In fact, it is straightforward
to modify the proof in the appendix to incorporate delayed responses.

The discussion of clinical applications has been the main focus of
this article because adaptive designs using urn models have received
much attention in statistics.  However, it is necessary to emphasize
that our results are very general and should also play an important
role in other areas as well.  For example, in quantum mechanics,
Niven and Grendar (2009) use the P\'{o}lya urn to understand the
generalized probability distribution for Maxwell-Boltzmann,
Bose-Einstein, and Fermi-Dirac statistics.  With different colors in
the urn, a ball is sampled, recorded and returned to the urn. Then,
$c$ balls of the same color are  added to the urn. In their
formulation, the choices of $c$ are $c > 0$, $c = 0$ and $c < 0$. As
$c < 0$ implies a decrease of the number of balls in the urn, it
would be interesting to explore the possibility of using the IMU
framework to avoid the distinction of balls of a particular type.

\vglue 0.12in

\section*{Appendix. Proofs}
\setcounter{equation}{0}
\renewcommand{\thesection}{A}

The outline of the proofs is as follows. First, we prove Theorem
\ref{th3}, which is our main result, and then Theorem \ref{th4}.
Finally we give a sketch of the proof of Theorem \ref{theoremGPU}.

 Recall that $\bm Z_{m-1}=(Z_{m-1,0}$, $Z_{m-1,1},$ $\ldots, Z_{m-1,K})$
 represent
the numbers of balls when the $m$-th subject arrives to be
randomized, $\bm Z_{m-1}^+=(Z_{m-1,0}^+,Z_{m-1,1}^+,$ $\ldots,$
$Z_{m-1,K}^+)$ are the non-negative numbers, and $|\bm
Z_{m-1}^+|=Z_{m-1,0}^+ + Z_{m-1,1}^+ +\ldots+ Z_{m-1,K}^+$.  Write
$\widetilde{\bm Z}_{m-1}=(Z_{m-1,1},\ldots,Z_{m-1,K})$. Because
every immigration ball is replaced, $Z_{m-1,0}^+=Z_{m-1,0}=Z_{0,0}$
for all $m$. Let $\bm X_m$ be the result of the $m$-th assignment,
where $X_{m,k}=1$ if the $m$-th subject is assigned to treatment $k$
and $0$ otherwise, $k=1,\ldots, K$. Then, $\bm N_n=(N_{n,1},\ldots,
N_{n,K})=\sum_{m=1}^n \bm X_m$. Further, we denote $\bm
a_m=(a_{m,1},\ldots,a_{m,K})$, and $\nu_m$ to be the number of draws
of type $0$ balls between the $(m-1)$-th assignment and the $m$-th
assignment.

Note that  between the $(m-1)$-th assignment and the $m$-th
assignment, we have drawn $\nu_m$ balls of type $0$. Accordingly, we
have added $a_{m-1,k}\nu_m$ balls of type $k$ to the urn. However,
when a ball of type $k$ is drawn, it is not replaced and
another $D_{m,kj}$ balls of type $j$ are added to the urn. Hence, the
change in the number of balls after the $m$-th assignment is
\begin{equation}
\label{eqPF1.3}
\widetilde{\bm Z}_m-\widetilde{\bm Z}_{m-1}=\bm a_{m-1} \nu_m +\bm X_m (\bm D_m-\bm I).
\end{equation}
 It follows that
 \vspace{-0.1in}
\begin{align} \label{eqPF1.2}
&\widetilde{\bm Z}_n-\widetilde{\bm Z}_0=\sum_{m=1}^n\bm a_{m-1} \nu_m +\sum_{m=1}^n\bm X_m
(\bm D_m-\bm I)\nonumber\\
=&\sum_{m=1}^n\bm a_{m-1} \nu_m -\bm N_n(\bm I-\bm H)+
\sum_{m=1}^n\bm X_m
(\bm D_m-\ep[\bm D_m])\nonumber\\
=&\bm a N_{n,0}+\sum_{m=1}^n(\bm a_{m-1}-\bm a) \nu_m
 -\bm N_n(\bm I -\bm H)+ \bm M_n,
\end{align}
where $N_{n,0}=\sum_{m=1}^n\nu_m$ is total number of draws of type
$0$ balls after  the $n$-th assignment, and
$\bm M_n =\sum_{m=1}^n \bm X_m (\bm D_m- \ep[\bm D_m])$ is a martingale.

\vglue 0.1in
 To prove  Theorem \ref{th3} we need two lemmas.  Their proofs will be given later.

\begin{lemma}\label{lem2} Suppose that the assumptions in Theorem \ref{th3} are satisfied.
 Then, for   $0<\delta_0<\frac{1}{2}-\frac{1}{2+\delta}$,
\begin{equation}\label{eqlem2.4}
Z_{n,k}=o(n^{1/2-\delta_0}) \quad a.s., \quad k=1,\ldots, K.
\end{equation}
\end{lemma}

\begin{lemma}\label{lem3} Suppose that the assumptions in Theorem \ref{th3} are satisfied. Then,
\begin{equation}\label{eqLILNast0}
N_{n,0}=n/s+O(\sqrt{n\log\log n})\;\; a.s.,
\end{equation}
\begin{equation}\label{eqLILNastk}
N_{n,k}=nv_k+O(\sqrt{n\log\log n}) \;\; a.s.,\;\;
k=1,\ldots, K,
\end{equation}
where $s=\bm a(\bm I-\bm H)^{-1}\bm 1^{\prime}$. Also, for each $k=1,\ldots,K$,
\begin{equation}\label{eqlem3.1}\widehat{\theta}_{n,k}\to \theta_k
\quad a.s.
\end{equation}
 and
\begin{equation}\label{eqapptheta}
 \widehat{\theta}_{n,k}-\theta_k=
\frac{Q_{n,k}}{nv_k}+o(n^{-1/2-\delta_0})\quad a.s.,
\end{equation}
where $Q_{n,k}=\sum_{m=1}^nX_{m,k}(\xi_{m,k}- \ep\xi_{m,k})$ is a
martingale and $\bm Q_n=(Q_{n,1},\ldots, Q_{n,K})$.
\end{lemma}

Now we begin the proof of Theorem \ref{th3}.  Consider the
$2K$-dimensional martingale $\{(\bm M_n, \bm Q_n),\Cal A_n; n\ge
1\}$, where $\Cal A_n=\sigma(\bm X_1,\ldots,\bm X_n, \bm
\xi_1,\ldots, \bm \xi_{n+1})$. According to (\ref{eqLILNastk}) we
have
\noindent
\begin{equation} \label{eqvarianceofM1}
\sum_{i=1}^n\ep[  (\Delta \bm M_i)^{\prime}\Delta\bm M_i | \Cal
A_{i-1}]  =  \sum_{k=1}^K N_{n,k}\bm \Sigma_k
 = n \bm \Sigma_{11} +O(\sqrt{n\log\log n}) ~ a.s.,\\
\end{equation}
\begin{equation} \label{eqvarianceofM2}
\sum_{i=1}^n\ep[ ( \Delta \bm Q_i)^{\prime} \Delta \bm Q_i | \Cal
A_{i-1}]
  =\bm \Sigma_{22}diag(\bm N_n)= n \bm \Lambda_{22}+O(\sqrt{n\log\log n}) ~ a.s., \\
  \end{equation}
  \begin{equation} \label{eqvarianceofM3}
\sum_{i=1}^n\ep[  (\Delta \bm M_i)^{\prime} \Delta \bm Q_i | \Cal
A_{i-1}] = \bm \Sigma_{12} diag(\bm N_n) = n \bm\Lambda_{12}
+O(\sqrt{n\log\log n})~ a.s.
\end{equation}
 By Corollary 1.1 of Zhang (2004),
 we can define  the $2K$-dimensional Wiener processes
 $(\bm W(t), \bm B(t))$ with variance-covariance matrix $\bm\Lambda$
 such that for some $\epsilon>0$,
\begin{equation}\label{eqappMQ}
\bm M_n=\bm W(n) +o(n^{1/2-\epsilon})\;\; a.s., \;\; \bm Q_n =\bm
B(n)+o(n^{1/2-\epsilon})\;\; a.s.
\end{equation}
 Without loss of generality, we assume that
$\epsilon\le \delta_0$, where $\delta_0$ is defined as it is in
Lemma \ref{lem2}.
 Next, we need
to show that $(\bm W(t), \bm B(t))$ satisfies (\ref{eqth3.2}).
  Combining
(\ref{eqPF1.2}) and (\ref{eqlem2.4})
  yields
\begin{align} \label{eqLILofN}
\bm N_n (\bm I-\bm H)-\bm a N_{n,0}
=\bm M_n+\sum_{m=1}^n(\bm a_{m-1}-\bm a)\nu_m+o(n^{1/2-\delta_0})\;\;
a.s.
 \end{align}

  Recall that $\bm A=(\bm I-\bm H)^{-1}(\bm I-\bm 1^{\prime}\bm v)$,
  $\bm v=\bm a(\bm I-\bm H)^{-1}/(\bm a(\bm I-\bm H)^{-1}\bm 1^{\prime})$
  and note that $\bm N_n\bm 1^{\prime}=n$, $\bm a\bm A=s\bm v(\bm I-\bm 1^{\prime}\bm v)=\bm 0$.
 According to (\ref{eqLILofN}),
\begin{equation}\label{eqproofth3.10}
\bm N_n-n\bm v=
 \big(\bm M_n+\sum_{m=1}^n(\bm a_{m-1}-\bm a)\nu_m\big)\bm  A+o(n^{1/2-\delta_0})\;\; a.s.
 \end{equation}

 For $\bm a_m-\bm a$,  due to (\ref{eqapptheta}) and (\ref{eqappMQ}),
\begin{align}
 \bm a_m-\bm  a=&(\widehat{\bm\theta}_{m}-\bm\theta)\frac{\partial \bm  a(\bm\theta)}{\partial\bm\theta}
 +O\big(\|\widehat{\bm\theta}_{m}-\bm\theta\|^2\big)
 \nonumber\\
 =&
\frac{\bm Q_m}{m}diag\big(\frac{1}{\bm v} \big)\frac{\partial \bm  a(\bm\theta)}{\partial\bm\theta}
+o(m^{-1/2-\delta_0}) \label{eqappofAbyQ}\\
  =&
\frac{\bm B(m)}{m}diag\big(\frac{1}{\bm v}\big)\frac{\partial \bm
a(\bm\theta)}{\partial\bm\theta}+o(m^{-1/2-\epsilon}). \nonumber
 \end{align}
Note that immigration occurs only when a type $0$ ball is drawn. Let
$\tau_m$ be the total number of draws when the $m$-th type $0$ ball
is drawn. At that time, $\tau_m-m$ subjects have been assigned and
the $(\tau_m-m+1)$-th subject arrives to be randomized. Hence, we add
$a_{(\tau_m-m+1)-1,k}$ balls of type $k$ to the urn, $k=1,\ldots, K$.
It follows that
$$ \sum_{j=1}^na_{j-1,k}\cdot \nu_j=\sum_{m=1}^{N_{n,0}}a_{\tau_m-m,k}, $$
$$i.e., \quad \sum_{m=1}^n(\bm a_{m-1}-\bm a)\nu_m=\sum_{m=1}^{N_{n,0}}(\bm a_{\tau_m-m}-\bm a).
$$
It is easily seen that $\tau_m=\min\{n:N_{n,0}\ge m\}+m$.
 Due to (\ref{eqLILNast0}),
$$ \tau_m-m=\min\{n:N_{n,0}\ge m\}=sm+O(\sqrt{m\log\log m}) \;\; a.s. $$
It follows that
 \begin{align*}
 \bm a_{\tau_m-m}-\bm  a=&\frac{\bm B(\tau_m-m)}{\tau_m-m}diag\big(\frac{1}{\bm v}
 \big)\frac{\partial \bm  a(\bm\theta)}{\partial\bm\theta}+o(m^{-1/2-\epsilon}) \\
 =&\frac{\bm B(sm)}{sm}diag\big(\frac{1}{\bm v}\big)\frac{\partial \bm
a(\bm\theta)}{\partial\bm\theta}+o(m^{-1/2-\epsilon})\;\;a.s.
 \end{align*}
Using (\ref{eqLILNast0}), we conclude that
\begin{align}\label{eqproofth3.12}
&\sum_{m=1}^n(\bm a_{m-1}-\bm a)\nu_m
=\sum_{m=1}^{N_{n,0}}\left(\frac{\bm
B(sm)}{sm}diag\big(\frac{1}{\bm v} \big)\frac{\partial \bm
a(\bm\theta)}{\partial\bm\theta}+o(m^{-1/2-\epsilon})\right)
\nonumber\\
&=\int_0^{n/s}\frac{\bm
B(sx)}{sx}dx\;diag\big(\frac{1}{\bm v} \big)\frac{\partial \bm
a(\bm\theta)}{\partial\bm\theta}+o(n^{1/2-\epsilon})
\nonumber\\
&=\int_0^{n}\frac{\bm
B(x)}{x}dx\;diag\big(\frac{1}{\bm v} \big)\frac{1}{s}\frac{\partial \bm
a(\bm\theta)}{\partial\bm\theta}+o(n^{1/2-\epsilon})\;\; a.s.
\end{align}
However, it is easily checked that
\begin{equation} \label{eqproofth3.13} \frac{1}{s}\frac{\partial \bm
a(\bm\theta)}{\partial\bm\theta}\bm A = \frac{\partial \bm
v(\bm\theta)}{\partial\bm\theta}. \end{equation} Combining
(\ref{eqappMQ})-(\ref{eqproofth3.13}) the proof of (\ref{eqth3.2}) is
complete. $\Box$

\smallskip
Three more lemmas are needed before we prove  Lemmas \ref{lem2} and
\ref{lem3}.

\begin{lemma}\label{lem4} Under Assumption \ref{Asp2} and $Z_{0,0}>0$, we have
$$Z_{n,k}^-=O(1)\quad a.s., \;\; k=1,\ldots,K.
$$
\end{lemma}
{\bf Proof}. $~$ Note that $|\bm Z_m^+|\ge Z_{0,0}>0$ for all $m$ and so that the balls with negative numbers have
no chance of being drawn. In addition, at most $C+1$ balls of each treatment type have the chance of being removed only when a   ball of the same type is drawn because of the Assumption \ref{Asp2}.
 It follows that $Z_{n,k}\ge -C-1$.
 $\Box$

\begin{lemma}\label{lem5}
  Let $\Cal F_n=\sigma(\bm X_1,\ldots,\bm X_n,\bm Z_1,\ldots, \bm Z_n)$
be the history sigma field, and $A_m=\sum_{k=1}^K a_{m,k}$.
Suppose that Assumption \ref{Asp2} is satisfied. Then,
$\underline{A}:=\min_m A_m>0$  implies
\begin{equation}\label{eqlem5.1}  \ep[\nu_n^p|\Cal F_{n-1}]\le c_p
\Big(\big(\sum_{k=1}^K Z_{n-1,k}\big)^-/\underline{A}\Big)^{p+1}\frac{Z_{0,0}}{|\bm
Z_{n-1}^+|}\;\; a.s., \;\; \forall\; p\ge 1,
\end{equation}
where $c_p>0$ is a random variable that is a function of $Z_{0,0}$
and $\min_mA_m$. Particularly,
\begin{equation}\label{eqlem5.2}
\min_m A_m>0 \; \text{ implies }   \;\ep[\nu_n^p|\Cal F_{n-1}]=O(1)\quad a.s.
\end{equation}
\end{lemma}

 {\bf Proof}.   The event $\{\nu_n=l\}$
means that when the $n$-th subject is assigned, we have drawn $l+1$
balls continuously in which the first $l$ balls is of type $0$ and
the last one is not. Hence, $\pr(\nu_n=0|\Cal F_{n-1})=1-Z_{0,0}/|\bm
Z_{n-1}^+|$, and for $l=1,2,\ldots$, \vspace {-0.1in} \noindent
\begin{equation}\label{eqprooflem5.0} \pr\big(\nu_n=l|\Cal
F_{n-1}\big) =\frac{Z_{0,0}}{|\bm Z_{n-1}^+|}
 \prod_{j=1}^{l-1}\frac{Z_{0,0}}{|(\bm Z_{n-1}+ j\bm a_{n-1})^+|}
\cdot\Big(1-\frac{Z_{0,0}}{|(\bm Z_{n-1}+l\bm a_{n-1})^+|}\Big).
\vspace{-0.1in}
\end{equation}
Obviously, $
\pr\big(\nu_n=l|\Cal F_{n-1}\big) \le  Z_{0,0}/|\bm
Z_{n-1}^+|$, $l\ge 1$.
 Note that
$$|(\bm Z_{n-1}+\bm j\bm
a_{n-1})^+|=Z_{0,0}+\sum_{k=1}^K(Z_{n-1,k}+ja_{n-1,k})^+ \ge
Z_{0,0}+\sum_{k=1}^K Z_{n-1,k}+jA_{n-1}.
$$
It follows that $\underline{A}>0$ and $\sum_{k=1}^K Z_{n-1,k}\ge
-L\underline{A}$ imply for $l \ge L$,
\begin{align}\label{eqprooflem5.1}
\pr\big(\nu_n=l|\Cal F_{n-1}\big)\le \frac{Z_{0,0}}{|\bm
Z_{n-1}^+|}\prod_{j=L}^{l-1}\frac{Z_{0,0}}{Z_{0,0}+(j-L)\underline{A}}
\le c_0\frac{Z_{0,0}}{|\bm Z_{n-1}^+|}e^{-2(l-L)},
\end{align}
where $c_0>0$ depends only on $\underline{A}$ and $Z_{0,0}$. So
\begin{align*}
\ep[\nu_n^p|\Cal F_{n-1}]\le \sum_{l=1}^Ll^p \frac{Z_{0,0}}{|\bm
Z_{n-1}^+|}+\sum_{l=L+1}^{\infty}l^pc_0\frac{Z_{0,0}}{|\bm
Z_{n-1}^+|}e^{-2(l-L)}
\le c_p L^{p+1}\frac{Z_{0,0}}{|\bm
Z_{n-1}^+|}.
\end{align*}
Taking $L=\big[(\sum_{k=1}^K
Z_{n-1,k})^-/\underline{A}\big]+1$ completes the proof of (\ref{eqlem5.1}).  (\ref{eqlem5.2})
follows from (\ref{eqlem5.1}) and Lemma \ref{lem4}. $\Box$

\begin{lemma}\label{lem6} Suppose that Assumptions \ref{Asp0}-\ref{Asp2} are satisfied.  Then
\begin{equation} \label{eqboundofa}
\min_{m,k} a_{m,k}>0 \;\; \text{ and } \max_{m,k}a_{m,k}<\infty
\;\; a.s.
\end{equation}
\end{lemma}

{\bf Proof}. $~$ \quad By Lemma A.4 of Hu and Zhang (2004a), we have
\begin{equation}\label{eqprooftheta2}
N_{n,k}\to \infty \;\; \text{ implies } \widehat{\theta}_{n,k}\to
\theta_k\;\; a.s.,\;\; k=1,\ldots,K.
\end{equation}
Then,  $a_k(\bm y)>0$ for any $\bm y$ on $
\text{closure}\{\widehat{\bm\theta}_m;m=1,2,\ldots\}
=\bigotimes_{k=1}^K\{\theta_k,\widehat{\theta}_{m,k};m=1,2,\ldots\}.
$ By the continuity of $a_k(\cdot)$, (\ref{eqboundofa}) is satisfied.
 $\Box$

 {\bf Proof of Lemma
\ref{lem2}}. \quad   By Lemma \ref{lem6},
\begin{equation}\label{eqboundofAA}\underline{A}=:\min_m A_m>0\;\; \text{ and
}\;\;\overline{A}=:\max_m A_m<\infty
\end{equation}
Note that $\widetilde{\bm Z}_n\bm 1^{\prime}=\sum_{k=1}^KZ_{n-1,k}$. By
(\ref{eqlem5.1}) and Lemma \ref{lem4}, $\ep[\nu_n|\Cal F_{n-1}]\le C_0
Z_{0,0}/|\bm Z_{n-1}^+|$. So, according to (\ref{eqPF1.3}) or
(\ref{eqPF1.2}), we have
\begin{align}\label{eqprooflem2.2}
\widetilde{\bm Z}_n\bm 1^{\prime}=
&\widetilde{\bm Z}_{n-1}\bm 1^{\prime}+\nu_n A_{n-1}-\bm X_n(\bm I-\bm H)\bm 1^{\prime}+\Delta \bm M_n\bm 1^{\prime}
 \nonumber\\
\le  &\widetilde{\bm Z}_{n-1}\bm 1^{\prime}+A_{n-1}\ep[\nu_n|\Cal
F_{n-1}]-\underline{h}
+ A_{n-1}(\nu_n-\ep[\nu_n|\Cal F_{n-1}])+\Delta \bm M_n\bm 1^{\prime}\nonumber\\
 \le & \widetilde{\bm Z}_{n-1}\bm 1^{\prime}+C_0\overline{A}\frac{Z_{0,0}}{|\bm
Z_{n-1}^+|}-\underline{h}
   +\Delta U_n \nonumber\\
   \le & \widetilde{\bm Z}_{n-1}\bm 1^{\prime}+\Delta U_n-\underline{h}\; /2, \; \text{
if } \widetilde{\bm Z}_{n-1}\bm 1^{\prime}\ge 2C_0\overline{A}
Z_{0,0}/\underline{h},
\end{align}
where $\underline{h}=\min_k(1-\sum_{j=1}^Kh_{kj})>0$. Here,
$U_n=\sum_{m=1}^nA_{m-1}(\nu_m-\ep[\nu_m|\Cal F_{n-1}])+ \bm M_n\bm
1^{\prime}$ is a real martingale.
 Let  $S_n=\max\{1\le j\le n: \widetilde{\bm Z}_j\bm 1^{\prime}<2C_0\overline{A} Z_{0,0}/\underline{h}
\;\}$, where $\max(\emptyset)=0$.  Then, according to
(\ref{eqprooflem2.2}),
\begin{align}\label{eqbeforelem2.24}
&\widetilde{\bm Z}_n\bm 1^{\prime}\le \widetilde{\bm Z}_{n-1}\bm 1^{\prime}+\Delta U_n -\underline{h}\;/2
 \le \ldots \nonumber \\
\le& \widetilde{\bm Z}_{S_n}\bm 1^{\prime}+\Delta U_{S_n+1}+\ldots+\Delta
U_n-(n-S_n)\underline{h}\;/2
\nonumber \\
 \le& |\bm Z_0|\vee \big(2C_0\overline{A} Z_{0,0}/\underline{h}\big) + U_n-U_{S_n}-(n-S_n)\underline{h}\;/2.
\end{align}
 For the martingale
$\{U_n,\Cal F_n;n=1,2,\ldots\}$, we have
$$  \ep[|\Delta U_n|^{2+\delta}|\Cal F_{n-1}]\le
 C+C\max_jA_j^{2+\delta}=O(1)
$$
due to  Assumption \ref{Asp2} and (\ref{eqlem5.2}). Accordingly, we
can show that
\begin{equation}\label{eqLILofU}
 U_n =O\big(\sqrt{n\log\log n}\big)\;\; a.s.,
\end{equation}
\begin{equation}\label{eqincofU}
\max_{m\le \sqrt{n\log n}}|U_{n-[\sqrt{n\log
n}]+m} -U_{n-[\sqrt{n\log
n}]} |=o(n^{\frac{1}{2+\delta}}\log n )\;\;a.s.
\end{equation}
If $n-S_n\ge \sqrt{n\log n}$, then for $n$ large enough
\begin{align*}
U_n -U_{S_n} -(n-S_n)\underline{h}\;/2\le
O\big(\sqrt{n\log\log n}\big)- \underline{h}\sqrt{n\log n}/ 2<0
\end{align*}
due to (\ref{eqLILofU}). Note that $n\ge S_n$. If $n-S_n<\sqrt{n\log
n}$, then
\begin{align*}
 U_n -U_{S_n} -(n-S_n)\underline{h}\;/2
 \le & ~2\max_{m\le \sqrt{n\log
n}}|U_{n-[\sqrt{n\log n}]+m} -U_{n-[\sqrt{n\log
n}]} |\\
=& ~o(n^{\frac{1}{2+\delta}}\log n )\;\;a.s.
\end{align*}
by (\ref{eqincofU}). It follows that
$ \sum_{k=1}^K Z_{n,k}\le o(n^{1/2-\delta_0})$
a.s.
 due to (\ref{eqbeforelem2.24}).
However,  $Z_{n,k}^-=O(1)$ a.s. by Lemma   \ref{lem4}.
(\ref{eqlem2.4}) is proved. $\Box$

\smallskip {\bf Proof of Lemma \ref{lem3}}. Recall
$Q_{n,k}=\sum_{m=1}^n X_{m,k}(\xi_{m,k}-\theta_k)$ $k=1,\ldots, K$,
and both $\{ M_{n,k},\Cal A_n; n\ge 1\}$ and $\{ Q_{n,k},\Cal A_n;
n\ge 1\}$ are martingales. According to the law of the iterated
logarithm for martingales, we have
\begin{equation} \label{eqLILofM}
 M_{n,k}=O(\sqrt{n\log\log n}) \;\;\text{ and }\;\; Q_{n,k}=O(\sqrt{n\log\log n})
 \quad a.s.
 \end{equation}

However,   for each $k=1,\ldots, K$,
\begin{equation}\label{eqprooftheta1}\widehat{\theta}_{n,k}-\theta_k=
\frac{Q_{n,k}+O(1)}{N_{n,k} +c_2}\;\;
a.s.
\end{equation}

(\ref{eqLILofN}) remains  true by Lemmas \ref{lem2}. By (\ref{eqLILofN}) and  (\ref{eqLILofM}) we have
\begin{equation} \label{eqlowboundofN} \bm N_n(\bm I-\bm H)= \sum_{m=1}^n \bm a_{m-1}\nu_m +o(n)
\;\;a.s. \end{equation}
Note that all elements of the vector $\sum_{m=1}^n \bm a_{m-1}\nu_m$ are between $\underline{a} N_{n,0}$ and $\overline{a} N_{n,0}$,
 where $\underline{a}=\min_{m,k}a_{m,k}$ and $\overline{a}=\max_{m,k}a_{m,k}$. Hence, it is obvious that $\liminf_{n\to \infty} N_{n,0}/n>0$ a.s.,
  because otherwise the limit of $\bm N_n/n$ may be $\bm 0$ which contradicts to $\bm N_n\bm 1^{\prime}=n$. On the other hand, the $k$-th element of
  $\bm N_n(\bm I-\bm H)$ does no exceed $(1-h_{kk})N_{n,k}$. It follows that   $\liminf_{n\to \infty} N_{n,k} / n>0$ a.s. by (\ref{eqlowboundofN}), which, together with
(\ref{eqprooftheta1}) and (\ref{eqLILofM}), implies
$$ \widehat{\theta}_{n,k}-\theta_k=
O\Big(\frac{Q_{m,k}+O(1)}{n}\Big)=
O\Big(\sqrt{\frac{\log\log n}{n}}\Big)\to 0\;\; a.s.
$$
(\ref{eqlem3.1}) is proved and also
\begin{equation}\label{eqLILofak} a_{m,k}-a_k=a_k(\widehat{\bm\theta}_m)-a_k(\bm\theta)
=O(\|\widehat{\bm\theta}_m-\bm\theta\|)=O\big(\sqrt{ (\log\log
m)/{m}}\big) \quad a.s.
\end{equation}
Hence, by Theorem 2.18 of Hall and Heyde (1980) it is easy to check that $\sum_{m=1}^n(a_{m-1,k}-a_k)(\nu_m-\ep[\nu_m|\Cal F_{m-1}])=o(\sqrt{n})$ a.s.  It follows that
\begin{align}\label{eqLILofaku}
&\sum_{m=1}^n(a_{m-1,k}-a_k)\nu_m=\sum_{m=1}^n(a_{m-1,k}-a_k)\ep[\nu_m|\Cal F_{m-1}]+o(\sqrt{n}) \nonumber\\
=&\sum_{m=1}^n O\Big(\sqrt{\frac{\log\log
m}{m}}\Big)O(1)+o(\sqrt{n})=O(\sqrt{n\log\log n}) \;\; a.s.
\end{align}
by (\ref{eqlem5.2}) and (\ref{eqLILofak}).
  Combining (\ref{eqLILofN}),
(\ref{eqLILofM}), and (\ref{eqLILofaku}) yields
$$
\bm N_n- N_{n,0}\bm a(\bm I-\bm H)^{-1} =O( \sqrt{n\log\log n}) \quad a.s.,
$$
which, together with $\bm N_n \bm 1^{\prime}=n$, implies (\ref{eqLILNast0}) and (\ref{eqLILNastk}).
 Then, combining (\ref{eqLILNastk}), (\ref{eqLILofM}), and (\ref{eqprooftheta1}) yields
$$ \widehat{\theta}_{n,k}-\theta_k=
\frac{Q_{n,k}+O(1) }{nv_k
+O(\sqrt{n\log\log n})}=\frac{Q_{n,k}}{nv_k}
+o(n^{1/2-\delta_0}) \;\; a.s.
$$
 (\ref{eqapptheta}) is proved, and the proof of Theorem \ref{th3} is completed. $\Box$

\smallskip
{\bf Proof of Theorem \ref{th4}}. Note that Assumptions \ref{Asp0} and
\ref{Asp2} are satisfied, and $Z_{0,0}>0$,  $\bm H\bm 1^{\prime}=\bm
1^{\prime}$. Similar to (\ref{eqprooflem2.2}),
\begin{align*}
\widetilde{\bm Z}_n\bm 1^{\prime}= &\widetilde{\bm Z}_{n-1}\bm
1^{\prime}+\nu_n A_{n-1}+\Delta \bm
M_n\bm 1^{\prime}
  \le  \widetilde{\bm Z}_{n-1}\bm 1^{\prime}+C_0\overline{A}\frac{Z_{0,0}}{|\bm
Z_{n-1}^+|}
   +\Delta U_n \nonumber\\
   \le & \widetilde{\bm Z}_{n-1}\bm 1^{\prime}+C_0\overline{A}/
   \sqrt{n}+\Delta U_n,  \; \text{
if } \widetilde{\bm Z}_{n-1}\bm 1^{\prime}\ge Z_{0,0}\sqrt{n}.
\end{align*}
It follows that
 \begin{align*}
\widetilde{\bm Z}_n\bm 1^{\prime} \le& \widetilde{\bm Z}_{S_n}\bm
1^{\prime}+\Delta U_{S_n+1}+\ldots+\Delta
U_n+C_0\overline{A}(n-S_n)/\sqrt{n}
\nonumber \\
 \le&   2C_0\overline{A}\sqrt{n} + U_n-U_{S_n}
  \le 2C_0\overline{A}\sqrt{n} + 2\max_{m\le n}|U_m|,
\end{align*}
where $S_n=\max\{1\le j\le n:
\widetilde{\bm Z}_j\bm 1^{\prime}<
 Z_{0,0}/\sqrt{n}
\;\}$ and  $\max(\emptyset)=0$. Hence,
$$ \widetilde{\bm Z}_n= O(\sqrt{n\log\log n})\; a.s. \;\; \text{ and } \;\;
=O_P(\sqrt{n}), $$
by the properties of a martingale and Lemma \ref{lem4}. So, by (\ref{eqPF1.2})  and the law of the iterated logarithm of martingales, it
follows that
\begin{align*}
\bm N_n (\bm I-\bm H) =&\bm M_n+\sum_{m=1}^n\bm
a_{m-1}\nu_m-\widetilde{\bm Z}_n+\widetilde{\bm Z}_0\sum_{m=1}^n\bm
a_{m-1}\nu_m+O(\sqrt{n\log\log n})\;\; a.s.
 \end{align*}
Multiplying by $\bm 1^{\prime}$ yields $\sum_{m=1}^n\nu_m A_{m-1}= O(\sqrt{n\log\log n})\;\;
a.s.$, and then $N_{n,0}= O(\sqrt{n\log\log n})$ a.s. and $\sum_{m=1}^n\bm a_{m-1}\nu_m = O(\sqrt{n\log\log n})$
a.s. by (\ref{eqboundofa}). So,
\begin{align*}
(\bm N_n -n\bm v)(\bm I-(\bm H-\bm 1^{\prime}\bm v))=\bm N_n(\bm I-\bm H)
=O(\sqrt{n\log\log n})\;\;
a.s.
 \end{align*}
It follows that $\bm N_n-n\bm v=O(\sqrt{n\log\log n})\;\; a.s.$
because $(\bm I-(\bm H-\bm 1^{\prime}\bm v))$ is invertible. The
proof of $\bm N_n-n\bm v=O_P(\sqrt{n})$ is similar. $\Box$

\smallskip

{\bf Proof of Theorem \ref{theoremGPU}.}  Recall (\ref{eqPF1.2}); we have
\begin{align} \widetilde{\bm Z}_n-\widetilde{\bm Z}_0=\sum_{m=1}^n\bm a_{m-1} \nu_m
 +\bm N_n(\bm H -\bm I)+ \bm M_n.
\end{align}
It follows
that
$|\widetilde{\bm Z}_n|=\widetilde{\bm Z}_n\bm 1^{\prime}\ge (\gamma-1)n +\bm M_n\bm 1^{\prime} $ by noticing  $\bm H\bm 1^{\prime}=\gamma\bm 1^{\prime}$.
Hence,
$$ \liminf_{n\to \infty} \frac{|\widetilde{\bm Z}_n^+|}{n}\ge \liminf_{n\to \infty} \frac{|\widetilde{\bm Z}_n|}{n}\ge \gamma-1>0\;\; a.s. $$
Without loss of generality we can thus assume that $|\widetilde{\bm Z}_n^+|\ge cn>0$ for all $n$. Then, the conclusion of Lemma \ref{lem4} remains true.
 By Lemma \ref{lem4},
$\widetilde{\bm Z}_m=\widetilde{\bm Z}_m^++O(1)$ a.s..
On the other hand, by (\ref{eqprooflem5.0}) we have
$$\pr(\nu_m=1|\Cal F_{m-1})=\frac{Z_{0,0}}{|Z_{m-1}^+|}\left(1-\frac{Z_{0,0}}{|(Z_{m-1}+\bm a_{m-1})^+|}\right)\le c/m\;\; a.s., $$
$$\pr(\nu_m\ge 2|\Cal F_{m-1})
=\frac{Z_{0,0}}{|Z_{m-1}^+|}\frac{Z_{0,0}}{|(Z_{m-1}+\bm a_{m-1})^+|}\le \Big(\frac{Z_{0,0}}{|Z_{m-1}^+|}\Big)^2\le c/m^2\;\;. $$
It follows that $\pr(\nu_m\ge 2\;\; i.o.)=0$ and
$\sum_{m=1}^n I\{\nu_m=1\}= O(\log^2 n) \;\; a.s.$
by Theorem 3.3.9 (ii) of Stout (1974). So by the assumption stated in
 Theorem \ref{theoremGPU}  that $0\le a_{m,k}\le C m^{1/2-\delta_0}$,
$$\sum_{m=1}^n\bm a_{m-1} \nu_m=O(\max_{m\le n} A_{m-1})\Big(\sum_{m=1}^n I\{\nu_m=1\}+O(1)\Big)= o(n^{1/2-\delta_0/2}) \;\; a.s. $$
which means that the  immigrated balls can be neglected.  In addition,
\begin{align*}
&\pr(X_{m,k}=1|\Cal F_{m-1})\\
=& \frac{\bm Z_{m-1,k}^+}{Z_{0,0}+|\widetilde{\bm Z}_{m-1}^+|}\left(1- \frac{\bm Z_{0,0}^+}{Z_{0,0}+|\widetilde{\bm Z}_{m-1}^+|}\right)
+\pr(X_{m,k}=1,\nu_m\ge 1|\Cal F_{m-1})\\
=& \frac{\bm Z_{m-1,k}^+}{ |\widetilde{\bm Z}_{m-1}^+|}+O\Big(\frac{1}{m}\Big)\;\; a.s.
\end{align*}
 It follows that
\begin{align*}
 &\widetilde{\bm Z}_n^+= \widetilde{\bm Z}_n+O(1)=\bm N_n(\bm H -\bm I)+ \bm M_n+o(n^{1/2-\delta_0/2})\\
 =& \sum_{m=1}^n(\bm X_m-\ep[\bm X_m|\Cal F_{m-1}])(\bm H-\bm I)+\bm M_n\\
  &+
 \sum_{m=1}^{n}\ep[\bm X_m|\Cal F_{m-1}](\bm H-\bm I)+o(n^{1/2-\delta_0/2})\\
  =& \sum_{m=0}^{n-1}(\bm X_m-\ep[\bm X_m|\Cal F_{m-1}])(\bm H-\bm I)+\bm M_n\\
  &+\sum_{m=1}^n\Big[\frac{\widetilde{\bm Z}_{m-1}^+}{ |\widetilde{\bm Z}_{m-1}^+|}+O\Big(\frac{1}{m}\Big)\Big](\bm H-\bm I)+o(n^{1/2-\delta_0/2})\\
 =& (\gamma-1)n \bm v+ (\gamma-1)\sum_{m=1}^n(\bm X_m-\ep[\bm X_m|\Cal F_{m-1}])\widetilde{\bm H}+\bm M_n \\
  &+\sum_{m=0}^{n-1}\frac{\widetilde{\bm Z}_m^+}{ |\widetilde{\bm Z}_m^+|}(\gamma-1)\widetilde{\bm H}+o(n^{1/2-\delta_0/2})\;\; a.s.
 \end{align*}
The  expansion for  $\widetilde{\bm Z}_n^+$ is  similar to that for
$\bm Y_n$ in (6.2) of  Zhang and Hu (pp. 1421-1324, 2009). Hence, the
rest of the proof is omitted.

\smallskip

\begin{center} ACKNOWLEDGEMENTS \end{center}
\vspace{-0.1in} Special thanks go to the anonymous referee and the associate editor for their constructive comments,
which led to a much improved
version of the paper.

\bigskip

 \baselineskip 15pt
\begin{center}
REFERENCES
\end{center}

\vspace{-0.1in} \footnotesize

\begin{enumerate}

\item[{[1]}]{\sc Athreya, K. B.} and {\sc Karlin, S.} (1968).
Embedding of urn schemes into continuous time branching processes
and related limit theorems. {\em Ann. Math. Statist.} {\bf 39}
1801--1817.

\item[{[2]}] {\sc Athreya, K. B.} and {\sc Ney, P. E.} (1972). {\em Branching
Processes}. Spring-Verlag, Berlin.

\item[{[3]}] {\sc Bai, Z. D.} and {\sc  Hu, F.} (1999).
Asymptotic theorem for urn models with nonhomogeneous generating matrices. {\em Stochastic
Process Appl.}  {\bf 80} 87--101.

\item[{[4]}]{\sc Bai, Z. D.} and {\sc Hu, F.} (2005).
 Asymptotics in randomized urn models.
 {\em Ann. Appl. Probab.} {\bf 15} 914--940.

 \item[{[5]}]{\sc  Bai, Z. D.,   Hu, F.}  and  {\sc Rosenberger, W. F. } (2002).
 Asymptotic properties of adaptive
designs for clinical trials with delayed response. {\em Ann.
  Statist.} {\bf 30} 122--139.

\item[{[6]}]
{\sc Bai, Z. D.,  Hu, F.} and {\sc Shen, L.} (2002). An Adaptive
Design for Multi-Arm Clinical Trials. {\em Journal of Multivariate
Analysis,} {\bf 81}, 1-18.

\item[{[7]}]
 {\sc Bhattacharya}, R. (2008).
 Urn-based response adaptive procedures and optimality.
 {\em Drug Information Journal} {\bf 42} 441--448.

\item[{[8]}]
 {\sc Beggs}, A. W. (2005).
 On the convergence of reinforcement learning.
 {\em Journal of Economic Theory} {\bf 122} 1--36.

\item[{[9]}]
 {\sc Bena\"{i}m} M., {\sc Schreiber} S. J. and {\sc Tarr\`{e}s}, P. (2004).
 Generalized urn models of evolutionary processes.
 {\em Ann. Appl. Probab.} {\bf 14} 1455--1478.

\item[{[10]}]
{\sc Donnelly, P.} and {\sc Kurtz, T. G.} (1996). The asymptotic
behavior of an urn model arising in population genetics. {\em
Stochastic Process Appl.}  {\bf 64} 1--16.

\item[{[11]}]
{\sc Durham, S. D., Flournoy, N.} and {\sc Li, W.}
 (1998).  Sequential designs for maximizing
the probability of a favorable response. {\em Canadian Journal of
Statistics} {\bf 3} 479--495.

\item[{[12]}]{\sc Hall, P. \& Heyde, C. C.} (1980).
{\em Martingale Limit Theory and its Applications}.
Academic Press, London.

\item[{[13]}]
{\sc Higueras, I., Moler, J., Flo, F.} and {\sc San Miguel, M.}
 (2006), Central Limit theorems for generalized P\'{o}lya urn models. {\em J. Appl. Prob.} {\bf 43} 938--951.

\item[{[14]}]
{\sc Hoppe, F. M.} (1984). P\'{o}lya-like urns and the Ewens'
sampling formula.  {\em Journal of Mathematical Biology} {\bf 20}
91--94.

\item[{[15]}]
{\sc Hu, F.} and {\sc Rosenberger, W. F.}  (2003).  Optimality,
variability, power:  evaluating response-adaptive randomization
procedures for treatment comparisons.  {\em J. Amer. Statist. Assoc.}
{\bf 98} 671--678.

\item[{[16]}]
{\sc Hu, F.} and {\sc Rosenberger, W. F.}  (2006).  {\em The Theory
of Response-Adaptive Randomization in Clinical Trials}. John Wiley
and Sons. Wiley Series in Probability and Statistics.

\item[{[17]}]
{\sc Hu, F., Rosenberger, W. F.} and {\sc Zhang, L.-X.}  (2006).
Asymptotically best response-adaptive randomization procedures. {\em
J. Statist. Plann. Inf.}  {\bf 136} 1911--1922.

\item[{[18]}]
 {\sc Hu}, F. and {\sc Zhang}, L. X. (2004a).
 Asymptotic
properties of doubly adaptive biased coin designs for
multi-treatment clinical trials.
 {\em Ann. Statist.} {\bf 32} 268--301.

\item[{[19]}]{\sc Hu, F.} and {\sc Zhang, L.-X.}  (2004b). Asymptotic normality of urn
models for clinical trials with delayed response. {\em Bernoulli}
{\bf 10} 447--463.

\item[{[20]}]
{\sc Ivanova, A.}  (2003).  A play-the-winner type urn model with
reduced variability.  {\em Metrika} {\bf 58} 1--13.

\item[{[21]}]
{\sc Ivanova, A.}  (2006). Urn designs with immigration: Useful
connection with continuous time stochastic processes. {\em  J.
Statist. Plann. Inf.}  {\bf 136} 1836-1844.

\item[{[22]}] {\sc Ivanova, A.}  and {\sc Flournoy, N.} (2001). A birth and death
urn for ternary outcomes: stochastic processes applied to urn
models. In {\em Probability and Statistical Models with
Applications} (Charalambides, Ch. A., Koutras, M. V. and
Balakrishnan, N., Eds.). Chapman \& Hall/CRC, 583--600.

\item[{[23]}]
{\sc Ivanova, A., Rosenberger, W. F., Durham, S. D. } and {\sc
Flournoy, N.} (2000).  A birth and death urn for randomized clinical
trials: asymptotic methods.  {\em Sankhya B} {\bf 62} 104--118.

\item[{[24]}] {\sc Janson, S.} (2004). Functional limit theorems for multitype
branching processes and generalized P\'olya urns. {\em Stochastic
Process. Appl.}  {\bf 110} 177--245.

\item[{[25]}]
 {\sc Johnson}, N. L. and {\sc Kotz}, S. (1977).
  {\em Urn Models and Their Applications}. Wiley, New York.

 \item[{[26]}]
 {\sc Kotz}, S. and {\sc Balakrishnan}, N. (1997).
  Advances in urn models during the
past two decades, In {\em Advances in Combinatorial Methods and
Applications to Probability and Statistics} (Balakrishnan, N., Ed.)
Birkh\"{a}user, Boston.

\item[{[27]}] {\sc Knoblauch, K., Neitz, M.} and {\sc Neitz, J.} (2006).
An urn model of the development of L/M cone ratios in human and
macaque retinas. {\em Visual Neuroscience} {\bf 23} 387--394.

\item[{[28]}] {\sc Matthews, E. E., Cook, P. F., Terada, M.} and {\sc Aloia, M. S.} (2010).
Randomizing research participants: Promoting balance and concealment
in small samples. {\em Research in Nursing and Health} {\bf 33}
243--253.

\item[{[29]}] {\sc Milenkovic, O.} and {\sc Compton, K. J.} (2004).
Probabilistic transforms for combinatorial urn models. {\em
Combinatorics, Probability and Computing} {\bf 13} 645--675.

\item[{[30]}] {\sc Niven, R. K.} and {\sc Grendar, M} (2010).
Generalized classical, quantum and intermediate statistics and the P\'{o}lya urn model.
{\em Physics Letters A} {\bf 373} 621--626.

\item[{[31]}]
{\sc Rosenberger, W. F., Stallard, N., Ivanova, A., Harper, C. N.}
and {\sc Ricks, M. L.}  (2001).  Optimal adaptive designs for binary
response trials.  {\em Biometrics} {\bf 57}  909--913.

\item[{[32]}]
{\sc Smythe, R.T.} (1996). Central limit theorems for urn models.
{\em Stochastic Process. Appl.} {\bf 65} 115--137.

\item[{[33]}]{\sc Stout, W. F.} (1974). {\em Almost sure convergence}.
 Academic Press, New York.

  \item[{[34]}] {\sc  Tamura, R. N., Faries, D. E.,  Andersen, J. S.} and {\sc Heiligenstein, J. H.} (1994). A case study of an adaptive clinical trial
  in the treatment of out-patients with depressive disorder. {\em J. Amer. Statist. Assoc.} {\bf 89}
768--776.

 \item[{[35]}] {\sc  Tymofyeyev, Y., Rosenberger, W.F.} and {\sc Hu, F.} (2007). Implementing
optimal allocation in sequential binary response experiments.{\em J. Amer. Statist. Assoc.} {\bf 102}
224--234.

\item[{[36]}] {\sc Wei, L. J.}  (1979).  The generalized P\'olya's urn
design for sequential medical trials.  {\em Ann. Statist.} {\bf 7}
291--296.

\item[{[37]}]
{\sc Wei, L. J.} and {\sc Durham, S. D.}  (1978).  The randomized
play-the-winner rule in medical trials.  {\em J. Amer. Statist.
Assoc.} {\bf 73}, 840--843.


\item[{[38]}] {\sc Zhang, L.J.} and {\sc Rosenberger W. F.} (2006).
Response-adaptive randomization for clinical trials with continuous outcomes. {\em Biometrics} {\bf 62} 562¨C569.

\item[{[39]}] {\sc Zhang, L.-X.} (2004). Strong approximations of martingale vectors and their
applications in Markov-chain adaptive designs. {\em Acta Math.
Appl. Sinica, English Series} {\bf 20}(2) 337--352.

\item[{[40]}] {\sc Zhang, L.-X.} and {Hu, F.} (2009). The Gaussian approximation for multi-color generalized Friedman¡¯s urn model.
 {\em  Science in China Series A: Mathematics}
{\bf  52}(6) 1305¨C1326.

\item[{[41]}] {\sc Zhang, L.-X., Chan, W. S., Cheung, S. H.} and {\sc Hu, F.} (2007). A
generalized urn model for clinical trials with delayed responses.
{\em Statistica Sinica} {\bf 17} 387--409.

\item[{[42]}] {\sc Zhang, L.-X., Hu, F.} and {\sc Cheung, S. H.} (2006). Asymptotic theorems of
sequential estimation-adjusted urn models. {\em Ann. Appl. Probab.}
{\bf 16} 340--369.

\item[{[43]}] {\sc Zhu, H.} and {\sc Hu, F.} (2009). Implementing optimal allocation in sequential continuous response
experiments. {\em  J. Statist. Plann.
Inf.}  {\bf 139} 2420--2430.
\end{enumerate}

\end{document}